\label{key}\documentclass[journal]{IEEEtran}
\usepackage{cite}
\usepackage{amsmath,amssymb,amsfonts}
\usepackage{algorithmic}
\usepackage{graphicx}
\usepackage{color}
\usepackage{epstopdf}
\usepackage{textcomp}
\usepackage{multirow}
\usepackage{array}
\usepackage{ragged2e}
\usepackage{makecell}
\usepackage{booktabs}
\usepackage{subfigure}

\def\BibTeX{{\rm B\kern-.05em{\sc i\kern-.025em b}\kern-.08em
		T\kern-.1667em\lower.7ex\hbox{E}\kern-.125emX}}
\begin{document}
	
	%\title{Energy Storage Enhanced Grid-Forming \\ Hydrogen Electrolyzer}

\title{Enhanced Hydrogen Electrolyzer with Integrated Energy Storage to Provide Grid-Forming Services for Off-Grid ReP2H Application}
\author{Yiwei~Qiu,~\IEEEmembership{Member, IEEE}, Jiahao~Hu, Yi~Zhou,~\IEEEmembership{Member, IEEE}, Jie Zhu,~\IEEEmembership{Graduate Student Member, IEEE}, \\ Li~Jiang, Shi~Chen, \IEEEmembership{Member, IEEE}, and Buxiang~Zhou,~\IEEEmembership{Member, IEEE} 
	\thanks{This work was supported by the National Natural Science Foundation
		of China under Grant 52377116 and 52307126, and the Natural Science Foundation of Sichuan
		Province (2024NSFSC0870). (Corresponding author: Yi Zhou.) }
	\thanks{All authors are with the
		College of Electrical Engineering, Sichuan University, Chengdu 610065, China
		(e-mail: zhouyi2320@scu.edu.cn). }}
\maketitle

\begin{abstract}
	This article proposes an energy storage-enhanced hydrogen electrolyzer (ESEHE) to provide grid-forming (GFM) services for off-grid renewable power to hydrogen (ReP2H) systems. Unlike conventional ReP2H systems that use a centralized energy storage (ES) plant, the proposed topology directly connects batteries to the DC buses of electrolysis rectifiers.
	A tailored virtual synchronous machine (VSM) control framework enables the electrolyzer to autonomously provide real and reactive power support. A coordinated frequency-splitting energy extraction strategy is designed to exploit both the battery and the electrolysis stack's electrical double-layer (EDL) effect on different timescales, maximizing active power support while mitigating battery and stack degradation.
	An adaptive equalization control strategy is further developed to balance the battery state of charge (SOC) among multiple ESEHEs operating in parallel, which optimizes energy distribution and extends battery life.
	Real-time simulations on StarSim validate the proposed topology and control strategies. Techno-economic analysis shows that, compared with conventional off-grid ReP2H systems based on a centralized ES plant, the ESEHE improves overall energy efficiency by 0.23\% {\color{black}and reduces the initial total converter investment cost by roughly 6\%, mainly due to the elimination of bidirectional AC/DC conversion and its associated losses in the centralized ES plant}.
\end{abstract}


\section{Introduction}
\label{sec:introduction}

\IEEEPARstart{A}{midst} the green energy transition, many countries regard renewable power to hydrogen (ReP2H) as a key pathway to decarbonize sectors such as power, energy, and chemicals \cite{li2025redesigning,TowardsRenewable-DominatedEnergySystems}. Driven by techno-economic concerns and policy mandates, ReP2H systems are required to use local renewable resources while minimizing adverse effects on power grid stability \cite{Flexibilityassessmentandaggregation}. Off-grid ReP2H systems, free from grid interconnection constraints, can be deployed in remote areas, islands, or regions with weak grid access. By directly consuming local renewables, they benefit from low electricity costs. However, without grid support, maintaining power balance becomes challenging, and careful trade-offs between safety and economic efficiency are required in system design and operation.

Off-grid ReP2H systems typically adopt AC architectures \cite{zhu2024exploring}, where a grid-forming (GFM) energy storage (ES) plant operates along with the electrolyzer to maintain voltage/frequency (V/F) stability and energy balance \cite{Variableperiodsequencecontrolstrategy, Grid-FormingConverters}. Their design must balance renewable power utilization and conversion efficiency to gain competitiveness in the green hydrogen market \cite{Agreenhydrogeneconomy,Zeng2025Planning}. In practice, ES implements GFM control via either VSM or constant V/F \cite{OptimizationofopticalstorageVSGcontrolstrategy, ModellingImplementation}, while renewables and electrolyzers synchronize through phase-locked loops (PLLs) \cite{DesignofPowerSystemStabilizer}.

In engineering, most off-grid ReP2H systems are supported by a centralized ES plant \cite{NASTASI2023117293}, but this has drawbacks:
1) vulnerability to single-point failures: a fault in the ES plant can cause system-wide loss of V/F support, followed by a total loss of load  %cascading failures, and load shedding 
\cite{Reviewofenergystoragesystemtechnologies};
2) efficiency loss: bidirectional power flow through the ES station's AC/DC converters increases energy losses;
3) implementation constraints: a centralized ES plant requires dedicated infrastructure, which limits adaptability in harsh environments such as islands or plateaus.

Beyond ES-based GFM solutions, electrolyzers have been studied for grid support. Dozein \textit{et al.} \cite{VirtualInertiaResponseandFrequencyControl, Grid-FormingServicesFromHydrogenElectrolyzers} exploited the electrolysis stack's electrical double-layer (EDL) effect for fast V/F regulation. Huang \textit{et al.} \cite{AnalyticalModelingandControl} developed dynamic frequency regulation models for alkaline (ALK) electrolyzers to adapt to wind power. Torres \textit{et al.} \cite{10478586} proposed a hybrid electrolyzer system with decentralized power-sharing control, and Zhao \textit{et al.} \cite{Comparativestudyofbattery-based} examined cooperative control between offshore wind and electrolyzers in off-grid mode.

However, electrolyzers as GFM loads face key challenges:
1) the limited EDL capacitance provides V/F support for only seconds; 
2) electrochemical and process limits load range and ramp rates in large-signal regulation \cite{Extendedloadflexibilityof}; 
3) long-term stability and coordination with ES in off-grid settings are still underexplored.

To overcome these limitations, research has been shifting from single-source GFM control to multi-source coordination to address the limitations of standalone electrolyzers. Torres et al. \cite{TorresSupercapacitor} proposed a hybrid system combining ALK electrolyzers with supercapacitors to buffer rapid fluctuations while the electrolyzer manages slower dynamics. Saha et al. \cite{Sahabatteries} turns to combining diverse ES units to ensure grid compliance. 
	
Consequently, prior research highlights three critical requirements for off-grid ReP2H systems. These include autonomous GFM support, coordinated control between loads and ES, and optimal configuration for large-scale systems \cite{zhu2024exploring}. To fulfill these requirements, we propose a new concept of the ES-enhanced hydrogen electrolyzer (ESEHE), with the main contributions as follows:

1) Inspired by the concept of ES-enhanced wind farms \cite{XIAO2025110776}, we propose a topology of the ESEHE, each of which integrates a dedicated battery, for off-grid ReP2H applications. This decentralized ES design enables long-duration, high-capacity GFM operation and avoids the bidirectional conversion losses of centralized ES, improving efficiency.

\begin{figure*}
	\centering
	\includegraphics[width=0.62\linewidth]{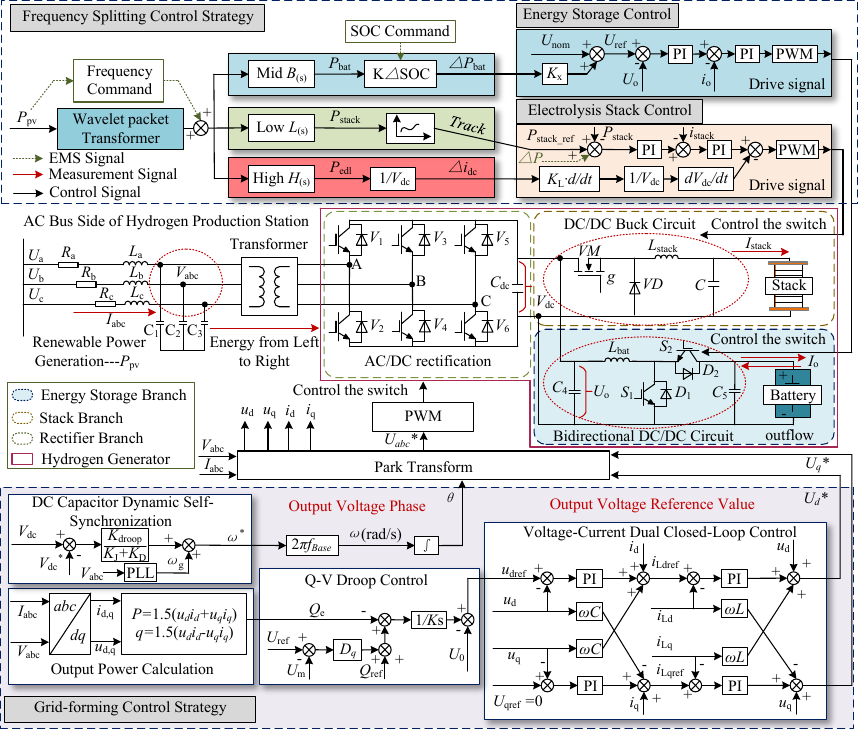}
	%\vspace{-0.2cm}
	\caption{Structure and control diagram of the proposed energy storage enhanced grid-forming hydrogen electrolyzer.}
	\vspace{-0.6cm}
	\label{fig:1}
\end{figure*}

2) {\color{black}A customized VSM controller coordinates energy between the ES and the electrolysis stack. Integrating this method directly into existing distributed and local controllers eliminates the need for extra infrastructure. This approach maximizes independent grid support while reducing investment and maintenance cost on standalone ES plants.}

3) The proposed ESEHE is validated via simulations of off-grid scenarios. Techno-economic analysis shows 0.23\% less energy loss {\color{black} and a 6\% reduction in total converter investment cost}  compared with conventional centralized ES-based systems. 

%{\color{black}guaranteeing superior reliability for remote deployments.}

The remainder of this paper is organized as follows. Section \ref{sec:Energy Storage Enhanced Grid-Forming Hydrogen Electrolyzer} presents the topology and controller design; Section \ref{sec:3} details the multi-timescale control, and Section \ref{sec:Case Analysis} gives real-time simulation results. Conclusions are drawn in Section \ref{sec:Conclusion}.

\section{Energy Storage Enhanced Grid-Forming Hydrogen Electrolyzer}
\label{sec:Energy Storage Enhanced Grid-Forming Hydrogen Electrolyzer}

%\subsection{Converter Control for Each Stage in the Electrolyzer}

Under VSM control, existing GFM electrolyzers \cite{VirtualInertiaResponseandFrequencyControl} use the electrolysis stack's EDL capacitance %and the DC-link capacitor
as the energy buffer.
However, this capability lasts only a few hundred milliseconds in proton exchange membrane (PEM) electrolyzers and a few seconds in ALK electrolyzers \cite{VirtualInertiaResponseandFrequencyControl}. Consequently, off-grid ReP2H systems cannot sustain VF stability during large disturbances without ES (see examples in Section \ref{subsec:Grid-Forming Capability Verification}).

To overcome this, we integrate a battery into the electrolysis rectifier's DC bus, as shown in Fig.~\ref{fig:1}. The energy buffering system, now comprising the DC-link capacitor, EDL capacitor, and battery, supports multiple timescales: the EDL capacitor responds to rapid transients, while the battery handles slower, larger fluctuations. 

To enable the GFM capability for the ESEHE, a customized VSM controller is designed for self-synchronization and real-time DC bus voltage regulation, stabilizing both the AC system and hydrogen production (see Section~\ref{sec:Energy Storage Enhanced Grid-Forming Hydrogen Electrolyzer}). Reliability and lifespan are enhanced via frequency-splitting control between the stack and ES, and SOC equalization among multiple electrolyzers (Sections \ref{subsec:Frequency-Splitting Power Control Method} and \ref{subsec:Equalization Control Strategy}). Direct battery integration also removes bidirectional conversion losses in centralized ES, improving efficiency (see Section~\ref{subsec:Techno-Economic Benefit Analysis}).

\subsection{Control Principles of the Proposed ESEHE}
\label{subsec:Control Principles of Energy Storage-Enhanced Grid-Forming Hydrogen Electrolyzer}

ReP2H systems require multi-timescale power control. At the hourly scale, an upper-level energy management system (EMS) schedules electrolyzer on/off states based on renewable forecasts \cite{zhu2024exploring, Extendedloadflexibilityof, SchedulingMultiple, buxiang2023optimal}. At the minute scale, the electrolytic load is adjusted to smooth medium- to low-frequency fluctuations. At the second scale, rapid VF stabilization is achieved using the EDL effect and DC bus capacitor. To integrate these multi-physics energy conversions in the ESEHE circuit, a customized GFM control framework is developed, as shown in Fig.~\ref{fig:1}.

%\label{subsec:Operating Modes of Various Converters}

Fig.~\ref{fig:1} illustrates the customized GFM control, integrating the ES battery, electrolysis stack, and DC bus capacitors to provide real/reactive power support for the AC bus while keeping the DC bus voltage stable. Power balance is maintained through precise reference tracking and adaptive GFM control. The subsystems are detailed below.

\subsubsection{AC/DC interface control}

Compared with conventional GFM control, the ESEHE's AC/DC interface is modified to emphasize DC bus voltage regulation. Reference \cite{AVirtualSynchronousControlforVoltage-Source} examines a GFM converter that self-synchronizes using the DC bus capacitor for voltage stabilization. Emulating generator rotor dynamics, the DC bus capacitor provides virtual inertia to the AC side, governed by:
\begin{equation}
	C_{\text{dc}} {dV_{\text{dc}}}/{dt} = I_{\text{dc}} - I_{\text{bat}},
	\label{1}
\end{equation}
where $C_{\text{dc}}$ is the DC capacitance; $V_\text{dc}$ is the DC bus voltage; $\textit{I}_\text{dc}$ is the current injected by the AD/DC converter; $\textit{I}_\text{bat}$ is the ES current. 
By controlling $\textit{I}_\text{dc}$, the following rotor-like dynamics can be emulated:
\begin{equation}
	J {d\omega}/{dt} = T_\text{m} - T_\text{e} - D \omega,
	\label{2}
\end{equation}
where $\textit{$J$}$ is the inertia constant; $\textit{D}$ is the damping coefficient; $\textit{T}_\text{m}$ is the mechanical torque; $\textit{T}_\text{e}$ is the electromagnetic torque. 

When a sudden change occurs in the system frequency, the DC capacitor $C_{\text{dc}}$ instantaneously regulates the DC bus voltage through (\ref{1}). The VSG converts the voltage deviation $\Delta V_\text{dc}$ into a frequency adjustment signal $\omega^{*}$ by emulating the generator rotor dynamics (\ref{2}), as follows:
\begin{equation}
	\omega^* = \omega_0 + K_{J} \int \Delta V_\text{dc} \, dt + K_\text{D} \Delta V_\text{dc},
	\label{3}
\end{equation}
where $\textit{K}_\textit{J}$ is the virtual inertia coefficient; $\textit{K}_\text{D}$ is the virtual damping coefficient; $\omega$$_\text{0}$ is the rated angular frequency. 

The AC/DC control, based on DC capacitor dynamics, is implemented in the \emph{Grid-Forming Control Strategy} block of Fig.~\ref{fig:1}, where active power reference adjusts in response to DC voltage deviation, inherently realizing power-voltage droop:
\begin{equation}
	P_{\text{ref}} = P_0 + K_{\text{droop}} \left( V_{\text{dc}}^* - V_{\text{dc}} \right),
	\label{4}
\end{equation}
where $\textit{P}_\text{ref}$ is the power reference; $\textit{P}_\text{0}$ is the initial setpoint, instructed by the upper-level EMS; $\textit{K}_\text{droop}$ is the droop factor; $\textit{V}$$_{\text{dc}}^{*}$ and $\textit{V}_\text{dc}$  are the rated and actual DC bus voltages. 

The frequency regulation signal is then obtained via dynamic capacitor feedback:
\begin{equation}
	\Delta \omega = \frac{1}{K_J s + K_\text{D}} (P_{\text{ref}} - P_\text{e}),
	\label{5}
\end{equation}
where $\Delta\omega$ is the frequency deviation; $\textit{P}_\text{e}$ is the rated power.

Eqs. \eqref{3}--\eqref{5} yield the frequency deviation $\Delta \omega = \frac{K_\text{droop}}{K_J s + K_\text{D}} (V_{\text{dc}}^*- V_\text{dc})$. This expression establishes the relationship between DC-link voltage and the motion of phase angle, which Fig. \ref{fig:1} depicts as the \emph{DC Capacitor Dynamic Self-Synchronous Control Strategy} block. The system generates the required reference voltage using a reactive power voltage droop loop and a dual closed-loop voltage current controller.

The phase angle $\theta$ is derived by integrating $\omega$. Combined with the voltage reference from the reactive power loop, this enables precise AC/DC converter regulation for autonomous GFM control.

\subsubsection{Energy storage control} 

The ES branch, equipped with a bidirectional DC/DC converter (see Fig. \ref{fig:1}), operates in either Buck or Boost mode depending on power availability.
When surplus power is available, the ES charges in Buck mode, where $\textit{S}_\text{1}$ is \textit{blocking} and $\textit{S}_\text{2}$ is \textit{conducting}.
The relation among the charging voltage $\textit{V}_\text{buck}$ and the DC bus voltage $\textit{V}_\text{dc}$ is:
\begin{equation}
	V_{\text{buck}} = D_{\text{buck}} \cdot V_{\text{dc}}(0 < D_{\text{buck}} < 1),
	\label{6}
\end{equation}
\noindent
where $\textit{D}_\text{buck}$ is the conduction duty cycle of switch $\textit{S}_\text{2}$. The corresponding charging power is:
\begin{equation}
	P_{\text{charge}} = V_{\text{buck}} \cdot I_{\text{buck}} = {V_{\text{dc}}^2 \cdot D_{\text{buck}}^2}/{\textit{R}_{\text{charge}}},
	\label{7}
\end{equation}
where $\textit{I}_\text{buck}$ denotes the charging current; and $\textit{R}_\text{charge}$ represents the equivalent resistance of the charging path. 

When power supply is insufficient, the ES can discharge in Boost mode to maintain hydrogen production. In this case, switch $\textit{S}_\text{2}$ is \textit{blocking} and $\textit{S}_\text{1}$ is \textit{conducting}, giving:
\begin{equation}
	V_{\text{dc}} = \frac{V_{\text{boost}}}{1 - D_{\text{boost}}},\  (0 < D_{\text{boost}} < 1),
	\label{8}
\end{equation}
where $\textit{D}_\text{boost}$ is the conduction duty cycle of switch $\textit{S}_\text{1}$. 
The corresponding discharge power is:
\begin{equation}
	\hspace{-3pt} P_{\text{discharge}} = V_{\text{boost}} \cdot I_{\text{boost}} = {V_{\text{dc}}^2 (1 - D_{\text{boost}}) D_{\text{boost}}}/{\textit{R}_{\text{discharge}}},
	\label{9} \hspace{-3pt} 
\end{equation}
where $\textit{I}_\text{boost}$ is the discharging current; and $\textit{R}_\text{discharge}$ denotes the equivalent resistance of the discharging path. The detailed coordination control described in Section \ref{subsec:Frequency-Splitting Power Control Method}.

\subsubsection{Electrolytic load control} 

\begin{figure}[!t]
	\centering
	\includegraphics[width=0.6\linewidth]{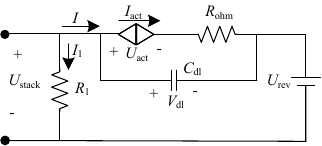}
	%\vspace{-0.2cm}
	\caption{Equivalent circuit model of the electrolysis stack.}
	\label{fig:2}
\end{figure}

The electrolysis stack's electrochemical characteristics can be modeled using an equivalent circuit, as shown in Fig. \ref{fig:2}.
	
The stack voltage consists of three parts: reversible overvoltage, activation overvoltage, and ohmic overvoltage \cite{Modelingstudyofefficiency}.
The reversible overvoltage $U_{\text{rev}}$, is determined by the Gibbs free energy change of the water-splitting reaction:
\begin{equation}
	U_{\text{rev}} = \frac{\Delta \textit{G}}{2\textit{F}} + \frac{\textit{RT}}{2\textit{F}} \ln \left( \frac{\textit{p}_{\text{H}_2} \sqrt{\textit{p}_{\text{O}_2}}}{a_{\text{H}_2\text{O}}} \right),
	\label{10}
\end{equation}
where $\textit{p}_{\text{H}_2}$ and $\textit{p}_{\text{O}_2}$ are the partial pressures of hydrogen and oxygen; $\textit{a}_{\text{H}_2\text{O}}$ is the activity of water; $\textit{R}$ is the gas constant; $\textit{T}$ is the stack temperature; and $\textit{F}$ is Faraday's constant. 

The activation overvoltage $U_{\text{act}}$, representing the voltage required to overcome the electrochemical activation barrier, can be approximated with a modified Tafel equation:
\begin{equation}
	U_{\text{act}} =\textit{v}_1\ln\left( {I}/{\textit{S}\textit{v}_2} + 1 \right),
	\label{11}
\end{equation}
where $\textit{v}_\text{1}$ and $\textit{v}_\text{2}$ are temperature-dependent parameters; $\textit{S}$ is the effective electrode area; and $\textit{I}$ is the stack current. 

The ohmic overvoltage $U_{\text{ohm}}$, representing voltage drop across internal resistance, is:
\begin{equation}
	U_{\text{ohm}} = I\textit{R}_{\text{ohm}}.
	\label{12}
\end{equation}

The total stack voltage is thus:
\begin{equation}
	U_{\text{stack}} = U_{\text{rev}} + U_{\text{act}} + U_{\text{ohm}}.
	\label{13}
\end{equation}
	
\begin{figure}[!t]
	\centering
	\includegraphics[width=0.75\linewidth]{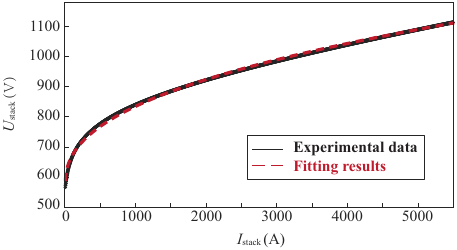}
	%\vspace{-0.2cm}
	\caption{U-I characteristic of the electrolysis stack.}
	\vspace{-0.4cm}
	\label{fig:d2}
\end{figure}

\begin{figure*}
	\centering
	\includegraphics[width=0.95\linewidth]{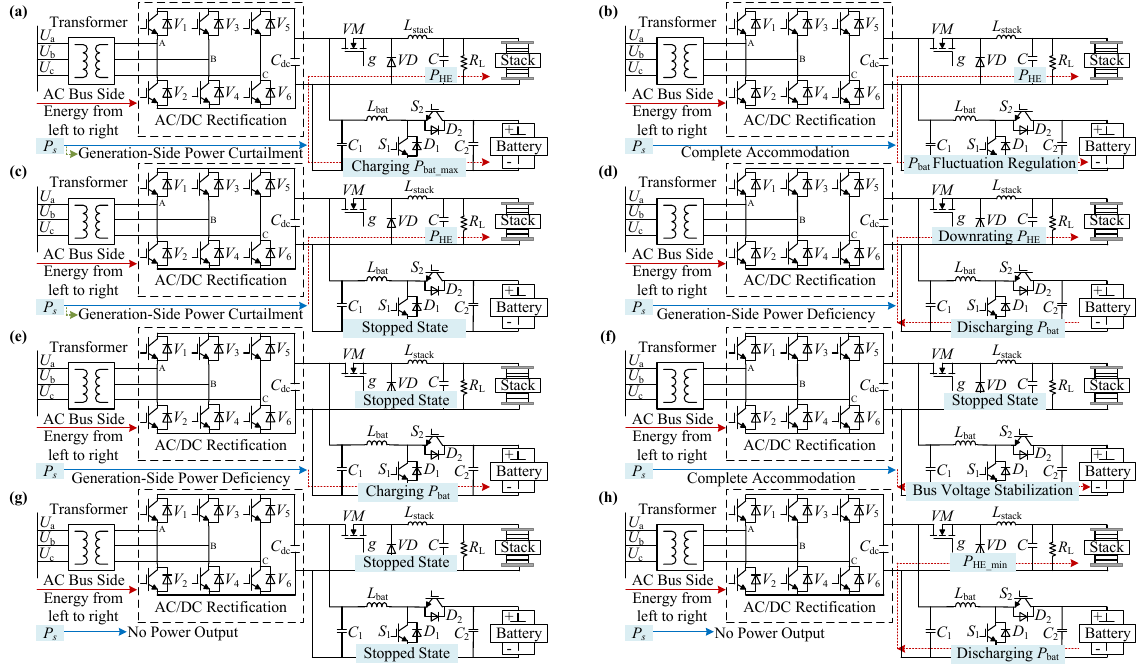}
	%\vspace{-0.2cm}
	\caption{Operational modes and energy flow diagram of the proposed ESEHE.}
	\label{fig:3}
	%\vspace{-0.5cm}
\end{figure*}

Eqs. \eqref{10}--\eqref{13} define the steady-state model of the electrolysis stack using the parameters listed in Table \ref{table:2} \cite{Modelingstudyofefficiency}, \cite{Power_controller_design}. Fig. \ref{fig:d2} illustrates the experimentally validated U-I characteristic derived from this model.

In addition to steady-state voltage, 
dynamic behavior from charge accumulation at the electrode-electrolyte interface, i.e., the EDL effect, must also be considered for fast regulation. The capacitance $C_\text{dl}$ depends on electrode material, surface morphology, and operating conditions, with dynamics:
\begin{equation}
	C_{\text{dl}} {dU_{\text{dl}}}/{dt} = I - I_{\text{act}},
	\label{14}
\end{equation}
\begin{equation}
	I_{\text{act}} = {U_{\text{act}}}/{\textit{R}_{\text{act}}},
	\label{15}
\end{equation}
\begin{equation}
	\textit{R}_{\text{act}} = {dU_{\text{act}}}/{dI},
	\label{16}
\end{equation}
where $\textit{I}_\text{act}$ is the current through the activation resistor; and $R_\text{act}$ represents its dynamically linearized value.Due to the complexity of measuring transients in large-scale systems, this study ensures reliability by adopting dynamic parameters from references \cite{Modelingstudyofefficiency}, \cite{Power_controller_design}, including experiment-verified EDL capacitance.

\begin{table*}[t]
	\centering
	\renewcommand{\arraystretch}{1.1} 
	\caption{State of Components under Different Operation Modes}
	%\vspace{-6pt}
	\label{table:1} 
	\begin{tabular}{@{}>{\centering\arraybackslash}m{0.13\columnwidth}  
			>{\centering\arraybackslash}m{0.4\columnwidth}>{\centering\arraybackslash}m{0.4\columnwidth}>{\centering\arraybackslash}m{0.5\columnwidth}>{\centering\arraybackslash}m{0.4\columnwidth}@{}} 
		\toprule 
		\toprule
		\textbf{Mode} & 
		\textbf{Power Supply} & 
		\textbf{Energy Storage} & 
		\textbf{Electrolysis Stack} & 
		\textbf{Power Relationship} \\ 
		\midrule 
		a & Partial power curtailment & Maximum-power charging & Maximum-power operation & $ P_{\mathrm{s}} > P_{\mathrm{HE}} $ \\
		b & Full consumption & Fluctuation regulation & Following renewable power output & $ P_{\mathrm{s}} = P_{\mathrm{HE}} + P_{\mathrm{bat}} $ \\
		c & Partial power curtailment & Shutdown & Rated-power operation & $ P_{\mathrm{s}} > P_{\mathrm{HE}} $ \\
		d & Insufficient output & Discharge & Reduced-power operation & $ P_{\mathrm{HE}} = P_{\mathrm{s}} + P_{\mathrm{bat}} $ \\
		e & Insufficient output & Charge & Shutdown & $ P_{\mathrm{s}} = P_{\mathrm{bat}} $ \\
		f & Full consumption & Charge and discharge & Shutdown & $ P_{\mathrm{s}} = P_{\mathrm{bat}} $ \\
		g & No output & Shutdown & Shutdown & $ P_{\mathrm{HE}} = 0 $ \\
		h & No output & Discharge & Minimum-power operation & $ P_{\mathrm{HE}} = \min (P_{\mathrm{e}}) $ \\
		\bottomrule 
		\bottomrule 
	\end{tabular}
	%\vspace{-0.4cm}
\end{table*}

In the \emph{Electrolysis Stack Branch} in Fig. \ref{fig:1}, the electrolytic power is regulated by the DC/DC converter's switch duty cycle $D$, which is as:
\begin{gather}
	\hspace{-6pt}\begin{cases} 
		\textit{L} {d\textit{i}_\text{L\_stack}}/{dt} = DV_{dc} - V_\text{stack} - \textit{i}_\text{L\_stack} \textit{R}_\text{L\_stack}, \vspace{4.5pt}\\
		C  dV_\text{stack}/{dt} = \textit{i}_\text{L\_stack} - I_\text{stack} = \textit{i}_\text{L\_stack} - {V_\text{stack}}/{\textit{R}_\text{act}},
	\end{cases} 	\hspace{-6pt}
	\label{17}
\end{gather}
where $\textit{L}_\text{stack}$ is the inductance; $\textit{i}_\text{L\_stack}$ and $\textit{I}_\text{stack}$ are the inductor and stack currents, respectively; $D$ is the duty cycle; $\textit{R}_\text{L\_stack}$ is the inductor's parasitic resistance; $C$ is the output capacitance.

The EDL effect can cause a response delay up to tens of seconds \cite{Power_controller_design}. To address this, a dual-loop control strategy is used. The outer loop generates the current reference $\textit{I}_\text{ref\_stack}$ from the power reference $\textit{P}_\text{ref\_stack}$:
\begin{equation}
	I_{\text{ref\_stack}} = {P_{\text{ref\_stack}}}/{U_{\text{stack}}};
	\label{18}
\end{equation}
the inner loop tracks $\textit{I}_\text{ref\_stack}$ with PI control to produce $D$:
\begin{equation}
	D = K_\text{p}\left( I_{\text{ref\_stack}} - I_{\text{stack}} \right) + K_\textit{i}\int\left( I_{\text{ref\_stack}} - I_{\text{stack}} \right) dt,
	\label{19}
\end{equation}
where $\textit{K}_\text{p}$ and $\textit{K}_\textit{i}$ are the proportional and integral gains of the current loop in electrolytic load control in Fig. \ref{fig:1}, related to the EDL dynamics and the charge/discharge rate limit of the ES. These gains are tuned via frequency-domain analysis and time-domain validation. 

To alleviate electrode degradation, the current slew rate is limited to $\leq$ 0.1 A/$\mu$s. 
With a bandwidth limit of $f_\text{c} \leq f_\text{sw}$/10, $\textit{K}_\text{p}$ and $\textit{K}_\textit{i}$ are optimized through Bode plot analysis to achieve a phase margin $>$65\textdegree, and gain margin $>$10 dB.

To actively release EDL energy for fast active power support, $D$ can be increased rapidly to boost $\textit{i}_\text{L\_stack}$, causing the EDL capacitor to discharge and provide a fast power boost.

Additionally, shunt current exists due to current leakage via lye channels \cite{QI2023233222}, represented by $\textit{R}_1$ and denoted as $\textit{i}_1$. Thus, the total stack current is:
\begin{equation}
	I_{\text{stack}} = \textit{i}_\text{L\_stack} + \textit{i}_1 + C_{\text{dl}} {dV_{\text{dl}}}/{dt}.
	\label{20}
\end{equation}

\subsection{Operating Modes}

The proposed ESEHE can switch among eight operation modes to adapt to different supply and load conditions, as shown in Fig. \ref{fig:3}, with details in Table \ref{table:1}\cite{jianlinS202411007}, explained as:

\emph{Modes a and b}: When $P_\text{s} > P_\text{HE}$, the surplus energy is stored in the battery; see Figs. \ref{fig:3}(a)--(b).

\emph{Mode c}: The ES is idle, the stack runs at rated power, and excess renewable energy is curtailed; see Fig. \ref{fig:3}(c).

\emph{Mode d}: When $P_\text{s}$ is insufficient for demanded hydrogen production, the ES discharges to fill the gap; see Fig. \ref{fig:3}(d).

\emph{Modes e and f}: The electrolysis stack is off, and all power exchange is via the ES branch ($P_\text{s} = P_\text{bat}$), functioning like a standalone ES station; see Figs. \ref{fig:3}(e)--(f).

\emph{Modes g and h}: With zero renewable generation or grid disconnection, the stack is either off (\emph{Mode g}) or runs at minimum power from the ES (\emph{Mode h}); see  Figs. \ref{fig:3}(g)--(h).

\section{Frequency-Splitting Control and SOC Equalization Strategy}
\label{sec:3}

\subsection{Frequency-Splitting Power Control}
\label{subsec:Frequency-Splitting Power Control Method}

Due to meteorological variability, off-grid ReP2H systems often face source-load mismatches. Conventional two-stage electrolysis rectifiers cannot absorb high-frequency and large-signal disturbances, and their controllers typically overlook the EDL effect, causing delayed power response \cite{Power_controller_design}.

In contrast, the proposed structure coordinates front-end (AC/DC interface) and back-end (electrolysis stack and ES branches) control. The back-end converter stabilizes the DC bus voltage, while the ES and electrolytic load are coordinated to respond to power fluctuations. To achieve this, a frequency-splitting controller decomposes the generation-side power spectrum to suppress high-frequency disturbances and enhance dynamic performance.

Existing studies predominantly employ supercapacitors or battery energy storage systems to absorb high-frequency components, resulting in increased system power losses and converter costs \cite{TorresSupercapacitor, Sahabatteries}. In contrast, this paper introduces a tri-band power splitter: the EDL capacitor replaces the supercapacitor to compensate millisecond-scale fluctuations;
the ES addresses second-to-minute variations, reducing battery degradation; and the electrolytic load absorbs low-frequency components.
This concept is shown at the top of Fig.~\ref{fig:1}. A filter decomposes the fluctuating power $P_\text{s}$ into high-, medium-, and low-frequency components \cite{ControlandPowerBalancing}:
\begin{equation}
	\begin{cases} 
		P_{\mathrm{h}} = \omega_{\mathrm{h}} \left( P_{\mathrm{s}} \right) & (f > 10~\text{Hz}), \\
		P_{\mathrm{m}} = \omega_{\mathrm{m}} \left( P_{\mathrm{s}} \right) & (1~\text{Hz} < f \leq 10~\text{Hz}), \\
		P_{\mathrm{L}} = \omega_{\mathrm{L}} \left( P_{\mathrm{s}} \right) & (f \leq 1~\text{Hz}),
		\label{21}
	\end{cases}
\end{equation}
where $\textit{P}_\text{h}$, $\textit{P}_\text{m}$, and $\textit{P}_\text{L}$ are extracted components; %, serving as the power references for the energy branches; 
$\omega_{\mathrm{h}}$, $\omega_{\mathrm{m}}$, and $\omega_{\mathrm{L}}$ denote the corresponding cutoff angular frequencies.

A high-pass filter extracts $P_\text{h}$, dynamically compensated by the EDL capacitor:
\begin{equation}
	H(s) = \frac{s^2}{s^2 + 2\zeta_{\mathrm{h}} \omega_{\mathrm{h}} s + \omega_{\mathrm{h}}^2},
	\label{22}
\end{equation}
with $\omega_{\mathrm{h}} = 20 \pi$~rad/s, i.e., $f_{\mathrm{h}}=10$~Hz, and damping ratio  $\zeta_{\mathrm{h}}=0.707$ for Butterworth flatness. The parameters ensure that the EDL capacitor can handle rapid fluctuations without accelerating stack degradation, suppress low-frequency interference, and minimize phase distortion.

A band-pass filter extracts $P_\text{m}$ for the ES branch, which operates on second-scale dynamics (e.g., $\mathrm{d}P/\mathrm{d}t \leq 0.1$~p.u./s):
\begin{equation}
	B(s) = \frac{2 \zeta_{\mathrm{m}} \omega_{\mathrm{m}} s}{s^{2} + 2 \zeta_{\mathrm{m}} \omega_{\mathrm{m}} s + \omega_{\mathrm{m}}^{2}},
	\label{23}
\end{equation}
where $\omega_{\mathrm{m}} = 2 \pi$ rad/s (i.e., $f_\text{m} = 1$ Hz), and $\zeta_{\mathrm{m}} = 0.707$. The bandwidth is $\Delta f \approx 1.414$~Hz. 
The resulting bandwidth $\Delta f \approx 1.414$~Hz enables medium-frequency power regulation in the 0.1-10 Hz range, limiting electrolytic load ramp rates to $\leq 0.1$~p.u./s and mitigating electrode and structural stress.

A low-pass filter isolates $P_\text{L}$ for the electrolytic reaction:
%(represented by $U_\text{rev}$, $U_\text{act}$, and $R_\text{ohm}$ in Fig. \ref{fig:2}), 
%ensuring relatively steady hydrogen production.
\begin{equation}
	L(s) = \frac{s^2}{s^2 + 2\zeta_{\mathrm{L}} \omega_{\mathrm{L}} s + \omega_{\mathrm{L}}^2},
	\label{24}
\end{equation}
with $\omega_{\mathrm{L}} = 0.2 \pi$~rad/s ($f_\text{L} = 1$~Hz) and $\zeta_{\mathrm{L}} = 0.770$. The cutoff matches the stack's load ramp capability, and the over-damped design avoids overshoot during start-up while maintaining passband gain within $\pm 0.05$ dB below $0.5$ Hz. 

The frequency band from 1 Hz to 10 Hz is partitioned based on the distinct physical response characteristics and operational constraints of the three coordinated energy components \cite{VirtualInertiaResponseandFrequencyControl},  \cite{Modelingstudyofefficiency}, \cite{fangbattery}. However, these cutoff thresholds may be dynamic rather than fixed. Considering the degradation and operational limits of different electrolyzer and energy storage types, determining optimal segmentation points warrants further study. Section \ref{subsec:Grid-Forming Capability Verification} presents the verification of this control strategy.

The frequency splitting control addresses degradation drivers across electrolyzer technologies. ALK units are sensitive to startup and shutdown cycles. A low-pass filter below 1 Hz aligns power commands with ramp limits to prevent costly shutdowns \cite{Flexibilityassessmentandaggregation, catalLIU}. Conversely, proton exchange membrane units are vulnerable to fluctuations. The system protects them by diverting components above 1 Hz to the ES and EDL effect, respectively \cite{Power_controller_design, TOMIC}. Furthermore, the adaptive SOC balancing strategy in Section \ref{subsec:Equalization Control Strategy} supports distributed management of the embedded ES units for large-scale systems.

Following \eqref{18}--\eqref{20}, the converter in the electrolysis stack branch regulates both outer-loop power and inner-loop current. The power loop tracks $P_\text{L}$:
\begin{align}
	I_{\text{ref\_stack}} &=  \frac{P_{\mathrm{L}}}{V_{\mathrm{stack}}} + K_{\mathrm{p}1} ( P_{\mathrm{L}} - V_{\mathrm{stack}} \text{\textit{i}}_{\mathrm{L\_stack}} ) \notag \\&+ K_{\textit{i}1} \int ( P_{\mathrm{L}} - V_{\mathrm{stack}} \text{\textit{i}}_{\mathrm{L\_stack}} )  dt,
	\label{25}
\end{align}
where $\textit{K}_\text{p1}$ and $\textit{K}_{\textit{i}\text{1}}$ are the PI gains  of the power loop in electrolysis stack control in Fig. \ref{fig:1}.

To suppress high-frequency disturbances in the EDL voltage $V_\text{dl}$ and reduce electrode stress, a feedforward term is added by linearizing the inductor equation:
\begin{gather}
	L {d\textit{i}_\text{L\_stack}}/{dt} \approx V_{\text{dc}} \Delta D - \Delta V_{\text{stack}},
	\label{26}\\
	\Delta D = K_{\text{ff}} \cdot {dV_{\text{dl}}}/{dt}.
	\label{27}
\end{gather}
Combining (\ref{20}), we have $\Delta V_\text{stack} \approx \Delta V_\text{dl}$ under high frequency. The feedforward gain is:
\begin{equation}
	K_{\text{ff}} = -{L_\text{stack}}/{V_{\text{dc}}}.
	\label{28}
\end{equation}

The modified duty cycle output is thus given as 
\begin{equation}
	D' = D + \Delta D.
	\label{29}
\end{equation}

The ES branch is connected in parallel to the electrolysis stack branch. Its dynamics are:
\begin{gather}
	\begin{cases}
		L_\text{b}  {d\textit{i}_\text{b} }/{dt} = D_\text{b} V_{\text{bat}} - V_\text{dc} - \textit{i}_\text{b} \textit{R}_\text{b}, \\
		\text{C}_\text{dc} {dV_\text{dc}}/{dt} = \textit{i}_\text{s}- \textit{i}_\text{L\_b} - \textit{i}_\text{b}  - {V_\text{dc}}/{\textit{R}_{\text{load}}},
	\end{cases}
	\label{30}
\end{gather}
where $\textit{L}_\text{b}$, $\textit{R}_\text{b}$, $\textit{i}_\text{b}$, $\textit{V}_\text{bat}$, $\textit{i}_\text{s}$, $\textit{i}_\text{L\_b}$, and $\textit{R}_\text{load}$ denote the inductance, resistance, ES current, terminal voltage, bus current, inductor current, and load resistance, respectively, of the ES branch.

The outer voltage loop stabilizes $V_\text{dc}$ and suppresses medium-frequency $P_\text{m}$ by setting $I_{\text{bat\_ref}}$:
\begin{equation}
	I_{\text{bat\_ref}} = K_{\text{p2}} \left( V_{\text{dc\_ref}} - V_{\text{dc}} \right) + K_{\text{\textit{i}\text{2}}} \int \left( V_{\text{dc\_ref}} - V_{\text{dc}} \right) dt,
	\label{31}
\end{equation}
where $V_{\text{dc\_ref}}$ is the voltage reference; and $K_\text{p2}$ and  $K_{\textit{i}2}$ are voltage-loop gains of the volatge loop in ES control in Fig. \ref{fig:1}. The control bandwidth is set above the disturbance range.

The inner current loop tracks $I_{\text{bat\_ref}}$ via duty cycle modulation, as follows:
\begin{equation}
	D_\text{b} = K_\text{p3} \left( I_\text{bat\_ref} - \textit{i}_\text{b} \right) + K_{\textit{i}\text{3}} \int \left( I_\text{bat\_ref} - \textit{i}_\text{b} \right) dt,
	\label{32}
\end{equation}
where $D_\text{b}$ is the ES branch duty cycle; and $K_\text{p3}$ and  $K_{\textit{i}3}$ are current-loop gains of the current loop in ES control in Fig. \ref{fig:1}. The current loop is tuned for fast response while limiting transient slew rates.

The frequency splitting control mitigates long-term degradation by restricting the electrolytic load to components below 1 Hz. Adherence to inherent ramp limits alleviates fatigue cracks and electrode damage \cite{catalLIU, huangiron}. While parameters such as ohmic overvoltage and capacitance drift over the stack lifespan, the proposed framework slows this process. As discussed in Section \ref{sec:Conclusion}, future work will integrate adaptive control schemes to compensate for these evolutionary changes.

\subsection{SOC Equalization Control for Multiple ESEHEs}
\label{subsec:Equalization Control Strategy}

\begin{figure}[!t]
	\centering
	\includegraphics[width=0.76\linewidth]{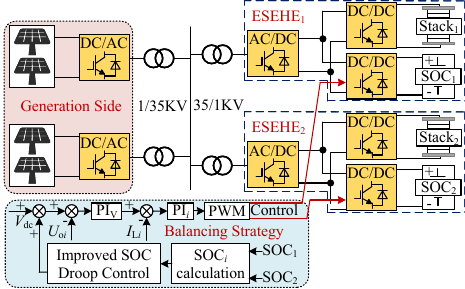}
	%\vspace{-7.5pt}
	\caption{ES balancing control strategy for multi-electrolyzer system.}
	\label{fig:4}
	%\vspace{-14pt}
\end{figure}

For large-scale applications, the limited capacity of a single stack (up to 7-10 MW) necessitates multiple electrolyzers operating in parallel. To achieve balanced utilization of the distributed ES units, an adaptive droop control is proposed, as shown in Fig. \ref{fig:4}. The goal is to equalize SOC across ES branches, thereby maximizing system-wide regulation capability and minimizing battery degradation. During charging, units with lower SOC are prioritized, with charging power reduced as SOC approaches the upper limit, and vice versa.

Specifically, inspired by \cite{DistributedOptimal}, an SOC-dependent function is superposed to the droop control's reference voltage. For the ES in the $i$th ESEHE, the voltage reference is given by:
\begin{equation}
	V_{\text{o}\textit{i}}^* = V_\text{dc} + f(\text{SOC}_\textit{i}),
	\label{33}
\end{equation}
where $f(\cdot)$ is a monotonically increasing function, designed as:
\begin{equation}
	f(\text{SOC}_\textit{i}) = \textit{k} \cdot \big( \text{e}^{\text{SOC}_\textit{i}^\textit{n}} - 1 \big) - \delta,
	\label{eq:placeholder}
\end{equation}
where $k$ and $n$ determine the equalization rate and non-linearity of the curve, while $\delta$ is a compensation factor to limit DC bus voltage deviation.

In parallel, an SOC-based adaptive droop coefficient is introduced to coordinate with the frequency-splitting control in Section \ref{subsec:Frequency-Splitting Power Control Method}. The adjustable ES capacity of each unit is defined by its SOC regulation range:
\begin{equation}
	\Delta \text{SOC}_\textit{i} = \text{SOC}_{\textit{i},\text{max}} - \text{SOC}_{\textit{i},\text{min}},
	\label{35}
\end{equation}
where $\text{SOC}_{\textit{i},\text{max}}$ and $\text{SOC}_{\textit{i},\text{min}}$ are the SOC limits of the ES unit; and $\Delta \text{SOC}_\textit{i}$ is the allowable SOC variation.

Based on $\Delta \text{SOC}_i$, the upward and downward droop coefficients in the frequency-power droop control are determined by the SOC status:
\begin{equation}
	\begin{cases} 
		\text{R}_{\text{up}} = \dfrac{(\text{SOC}_{\textit{i},\text{max}} - \text{SOC}_\textit{i})\Delta\text{SOC}_{\textit{i}}}{\sum\nolimits_{\textit{j} \in \Omega} ((\text{SOC}_{\textit{j},\text{max}} - \text{SOC}_\textit{j})\Delta\text{SOC}_{\textit{j}})}, \\
		\text{R}_{\text{down}} = \dfrac{(\text{SOC}_\textit{i} - \text{SOC}_{\textit{i},\text{min}})\Delta\text{SOC}_{\textit{i}}}{\sum\nolimits_{\textit{j} \in \Omega} ((\text{SOC}_\textit{j} - \text{SOC}_{\textit{j},\text{min}})\Delta\text{SOC}_{\textit{j}})},
	\end{cases}
	\label{36}
\end{equation}
where $\Omega$ is the set of parallel electrolyzers; and $\Delta \text{SOC}_i$ acts as a weighting factor. This design assigns larger droop coefficients for charging to units with lower SOC, and smaller ones for discharging to units with higher SOC.

\begin{figure*}[tb]
	\centering
	\subfigure[]{
		\begin{minipage}{.3\linewidth} 
			\centering
			\includegraphics[scale=0.84]{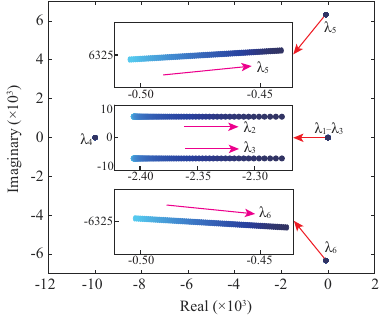}
			\label{d6:a}
		\end{minipage}
	}
	\subfigure[]{
		\begin{minipage}{.3\linewidth}
			\centering
			\includegraphics[scale=0.84]{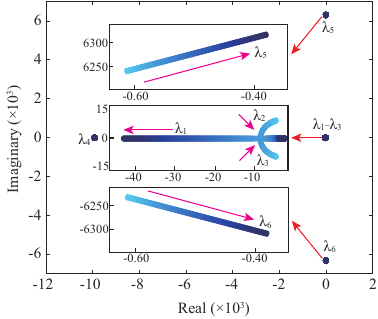}
			\label{d6:b}
		\end{minipage}
	}
	\subfigure[]{
		\begin{minipage}{.3\linewidth} 
			\centering
			\includegraphics[scale=0.84]{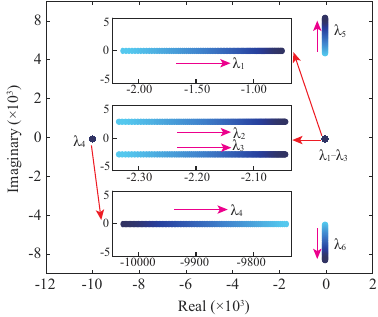}
			\label{d6:c}
		\end{minipage}
	}
	\caption{Small-signal stability verification via eigenvalue analysis of the proposed ESEHE system. (a) Root locus under varying load conditions. (b) Root locus under VSDC parameter variations. (c) Root locus under DC/DC parameter variations.}
	\label{fig:d6}
	%\vspace{-0.5cm} 	
\end{figure*}

It should be noted that this control requires an upper-level EMS to provide steady-state scheduling references (e.g., 96-point daily profiles) for the electrolytic load, the AC/DC active power, and the AC voltage. Existing studies on power management for multi-electrolyzer systems focus mainly on conventional HEs \cite{zhu2024exploring, Extendedloadflexibilityof, SchedulingMultiple, buxiang2023optimal}. Further research is needed for the long-term operation of the proposed ESEHEs.

\subsection{Small-Signal Stability Analysis}

We establish a small-signal state-space model to verify the stability of the coupling between VSM control and DC bus dynamics. This analysis isolates the rapid dynamics of the GFM loop and electrolysis stack. The ES branch is treated as a constant source due to its slower timescale operation, detailed in Section \ref{subsec:Frequency-Splitting Power Control Method}.

The non-linear system is linearized at the equilibrium. The state vector is defined as $\Delta x = [\Delta V_\text{dc}, \Delta \xi, \Delta \delta, \Delta I_\text{stack}, \Delta \gamma]^T$, which includes DC bus voltage deviation, the VSM inertia loop intermediate state, power angle, stack current, and the integral state of the stack current controller, respectively.

The linearized dynamics of the VSM control (DC voltage synchronization control, VSDC, introduced in Section \ref{subsec:Control Principles of Energy Storage-Enhanced Grid-Forming Hydrogen Electrolyzer}) and phase synchronization, derived from \eqref{3}--\eqref{5}, are given by:
\begin{equation} \begin{cases}
		{d\Delta\xi}/{dt} = \Delta V_\text{dc}, \\
		\Delta\omega = K_{J}\Delta\xi + K_\text{D}\Delta V_\text{dc}.
	\end{cases} 
	\label{37}
\end{equation}

The variation in AC active power is expressed as:
\begin{equation} 
	\Delta P_\text{ac} = K_{\delta}\Delta\delta, 
	\label{38}
\end{equation}
where $K_{\delta}$ is the synchronizing coefficient determined by the grid impedance and operating voltage.

The dynamics of the DC link and the electrolysis stack circuit, derived from \eqref{1} and \eqref{17}, are linearized as:
\begin{gather} \begin{cases}
		C_\text{dc}{d\Delta V_\text{dc}}/{dt} = \Delta I_\text{dc} - (D_{0}\Delta I_\text{stack} + I_\text{stack0}\Delta D), \\
		L_\text{stack}{d\Delta I_\text{stack}}/{dt} = D_{0}\Delta V_\text{dc} + V_\text{dc0}\Delta D - R_\text{stack}\Delta I_\text{stack},
	\end{cases}
	\label{39}
\end{gather}
where the subscript $0$ denotes steady-state values (duty cycle $D_0$, stack current $I_\text{stack0}$, and DC voltage $V_\text{dc0}$) at the equilibrium; $R_\text{stack}$ represents the equivalent resistance of the electrolysis stack branch linearized at the operating point; $\Delta I_\text{dc}$ is the small-signal current injection from the AC/DC rectifier.

The small-signal duty cycle $\Delta D$ is governed by the PI current controller:
\begin{gather} 
	\Delta D = K_\text{p}(I_\text{ref\text{\_}stack} - \Delta I_\text{stack}) + K_\text{i}\Delta \gamma, \label{40} \\
	{d\Delta \gamma}/{dt} = - \Delta I_\text{stack}. \label{41}
\end{gather}

By combining \eqref{37}--\eqref{41}, the small-signal state-space model is derived as $\dot{\Delta x} = A\Delta x + B\Delta u$. %, where $A$ is the system matrix.
Eigenvalue analysis of $A$ confirms stability. Fig. \ref{d6:a} presents the root locus as $I_\text{stack}$ varies from 1 kA to 5 kA, showing that dominant poles remain in the left half plane to ensure sufficient damping. Additionally, Fig. \ref{d6:b} demonstrates robustness to VSDC parameter changes with $K_\text{droop}$ ranging between 0.002 and 0.01, confirming stable GFM operation. Finally, Fig. \ref{d6:c} illustrates the impact of the DC/DC controller gain $K_\text{p3}$. These results indicate strong robustness where DC load dynamics do not destabilize the overall system.

\section{Case Study}
\label{sec:Case Analysis}

\subsection{Simulation Setup}
\label{subsec:Simulation Setup}

To validate the performance of the proposed ESEHE, we directly connect it to a 5 MW photovoltaic (PV) unit operating under conventional grid-following (GFL) control. The electrolysis stack is rated at 5 MW, with an AC/DC converter capacity of 6.63 MVA and an ES capacity of 0.8 MWh \cite{zhu2024exploring}. Key parameters are listed in Table \ref{table:2}.

\begin{figure}[!t]
	\centering
	\includegraphics[width=0.75\linewidth]{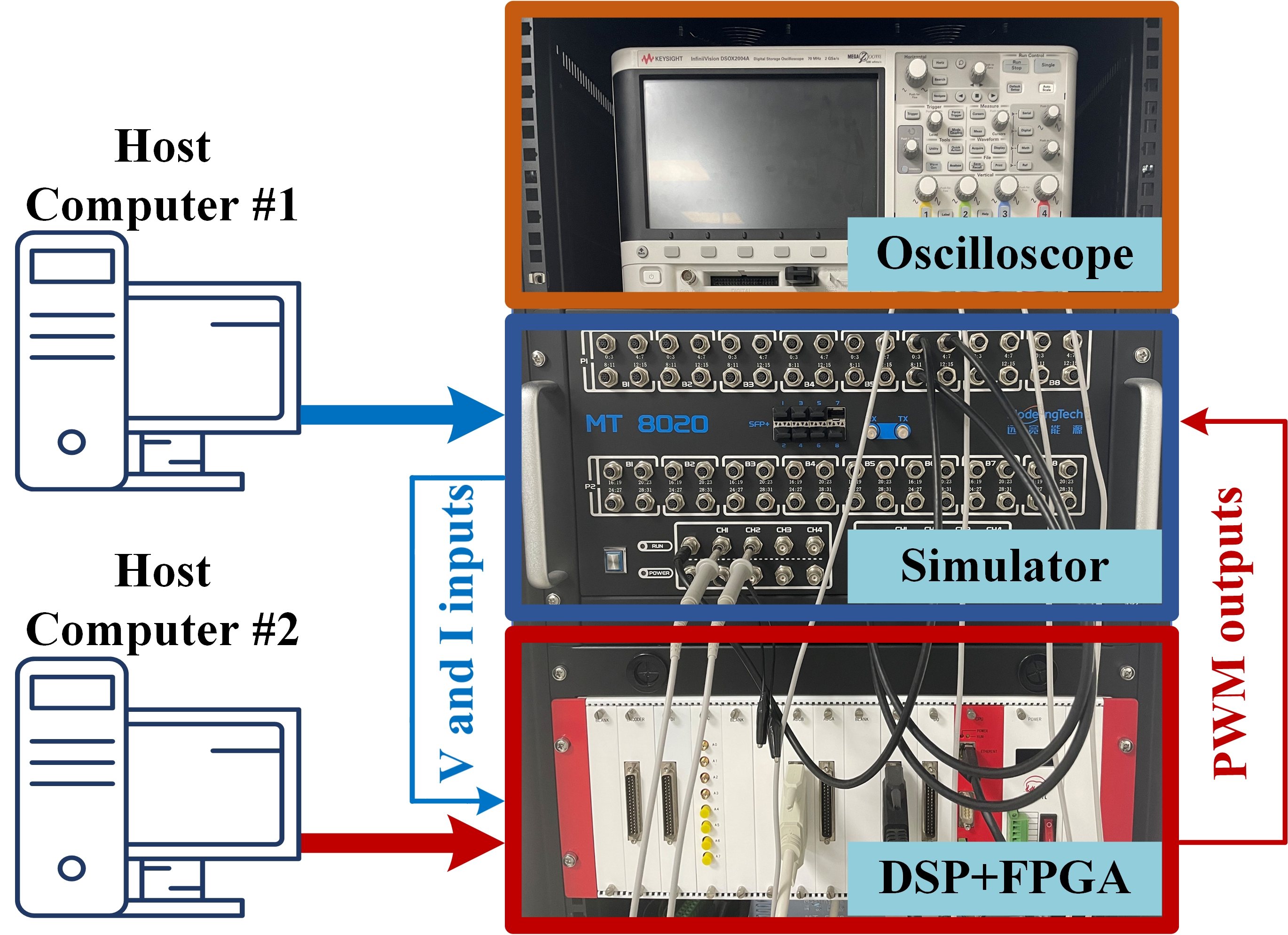}
	%\vspace{-4.5pt}
	\caption{Hardware-in-the-loop (HIL) platform used for simulation analysis.}
	\label{fig:5}
	%\vspace{-2pt}
\end{figure}

\begin{figure*}[!t]
	\centering
	\subfigure[]{
		\begin{minipage}[t]{.3\linewidth}
			\centering
			\includegraphics[scale=0.9]{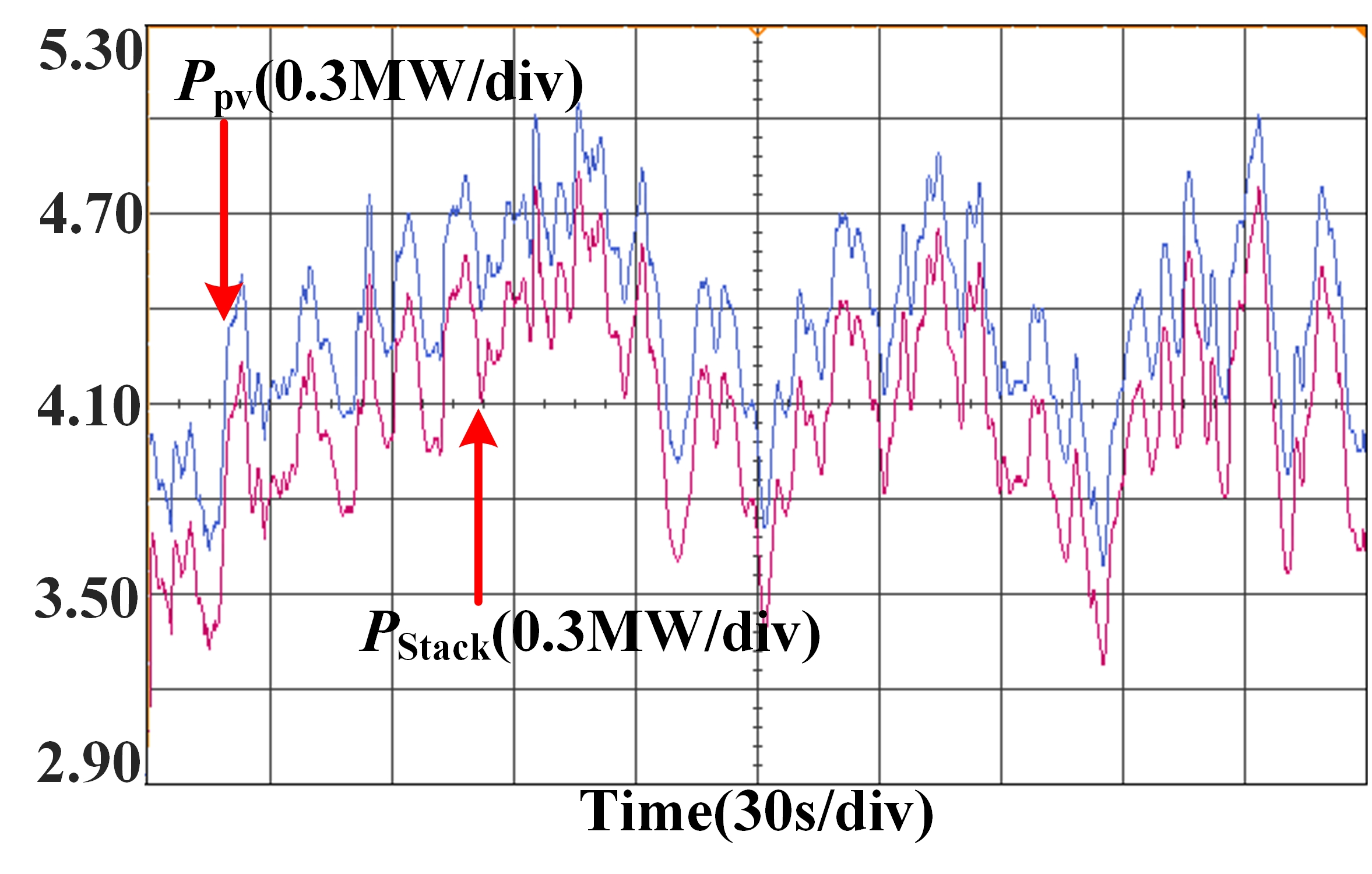}
			%\vspace{-1cm}
			\label{6:a}
		\end{minipage}
	}
	\subfigure[]{
		\begin{minipage}[t]{.3\linewidth}
			\centering
			\includegraphics[scale=0.9]{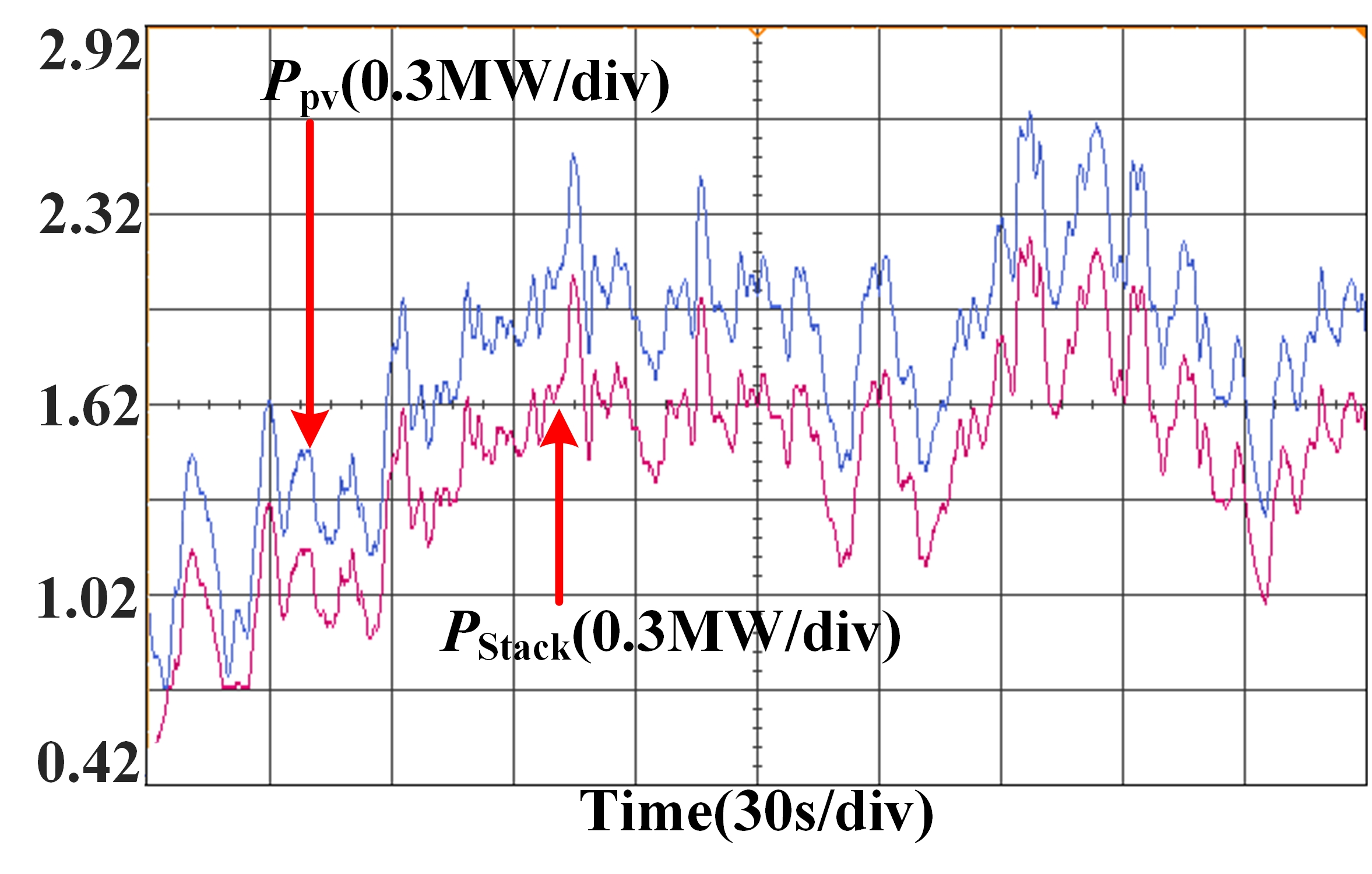}
			%\vspace{-1cm}
			\label{6:b}
		\end{minipage}
	}
	\subfigure[]{
		\begin{minipage}[t]{.3\linewidth}
			\centering
			\includegraphics[scale=0.9]{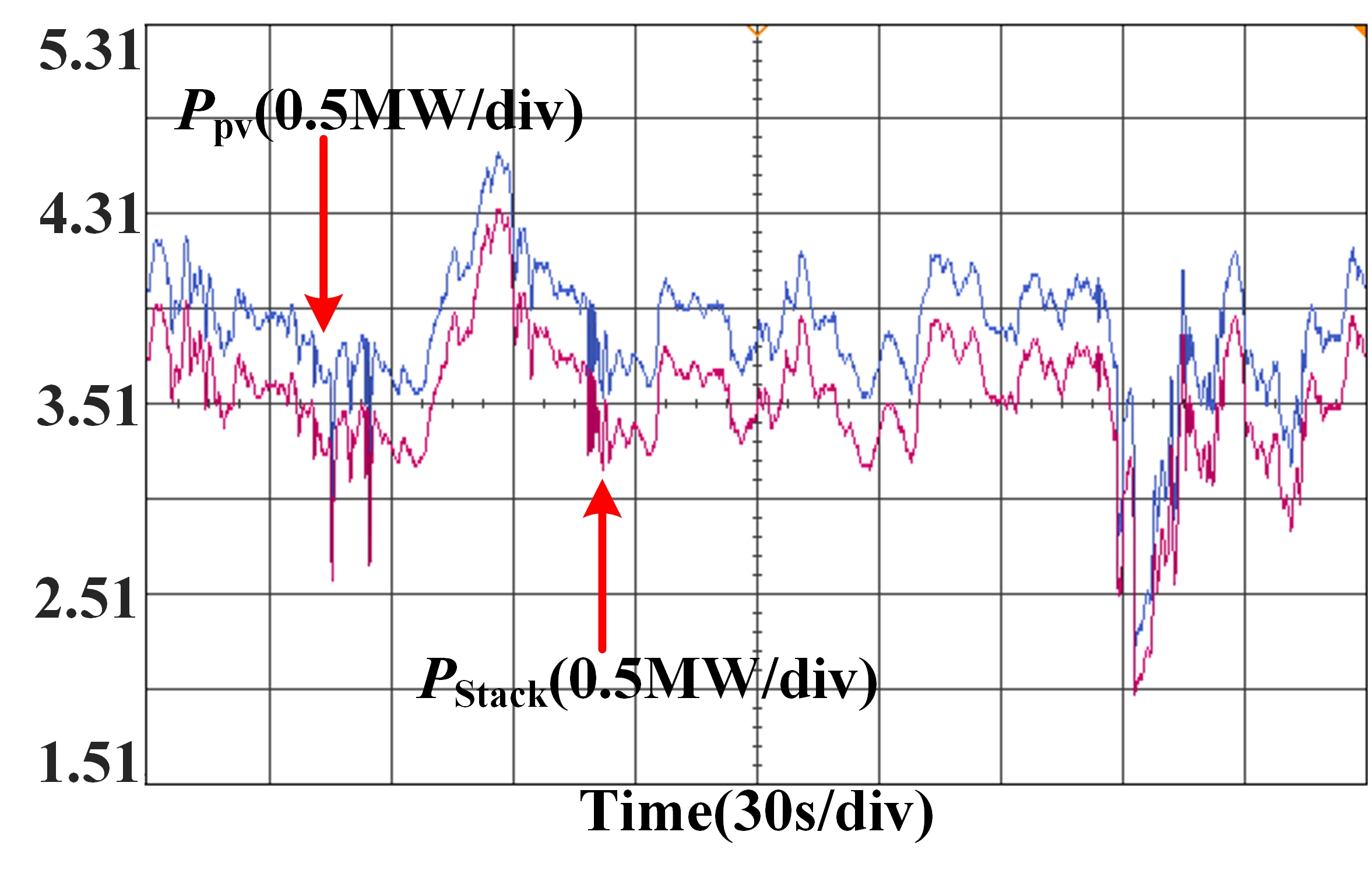}
			%\vspace{-1cm}
			\label{6:c}
		\end{minipage}
	}
	\subfigure[]{
		\begin{minipage}[t]{.3\linewidth}
			\centering
			\includegraphics[scale=0.9]{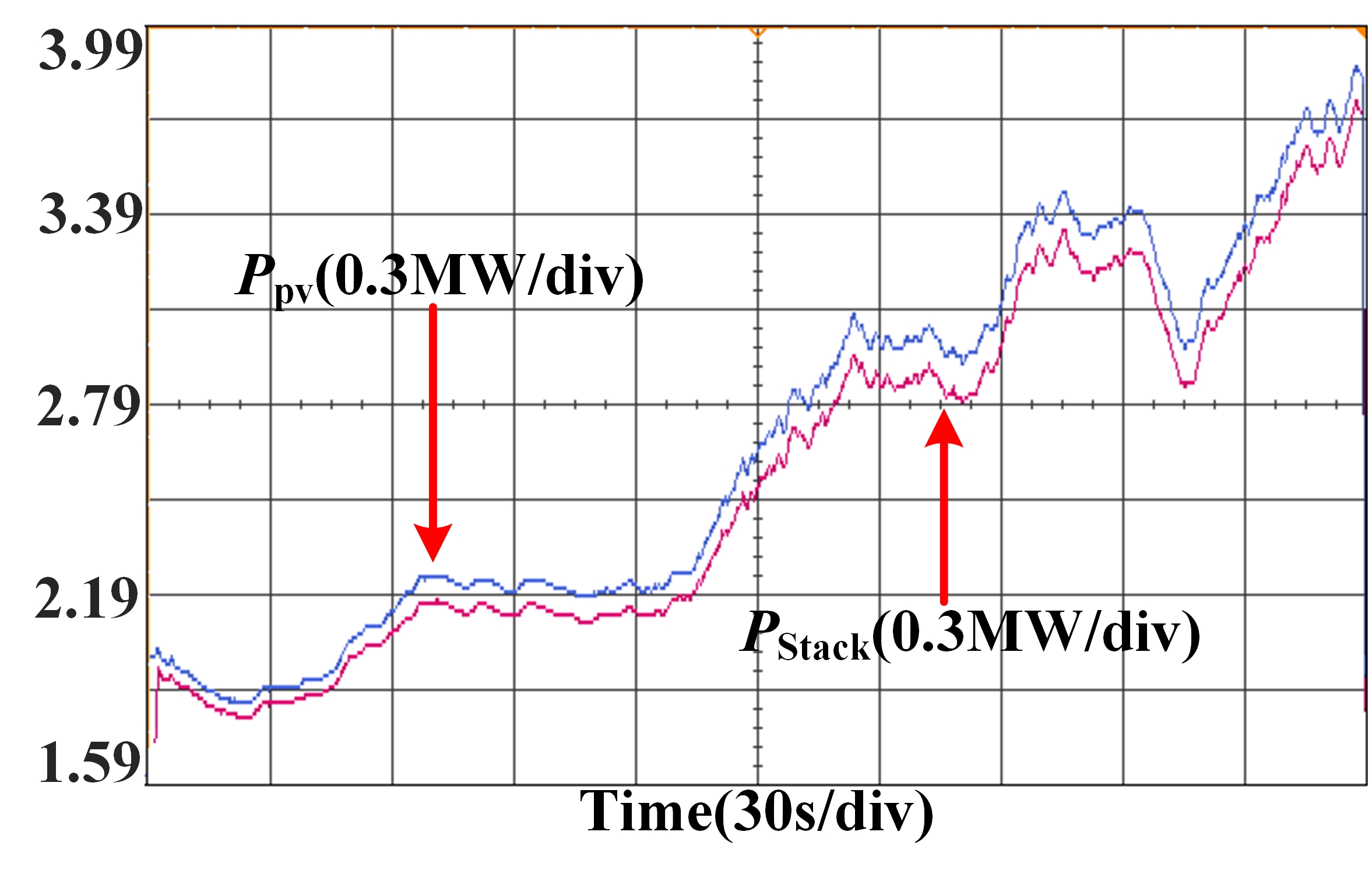}
			%\vspace{-1cm}
			\label{6:d}
		\end{minipage}
	}
	\subfigure[]{
		\begin{minipage}[t]{.3\linewidth} 
			\centering
			\includegraphics[scale=0.9]{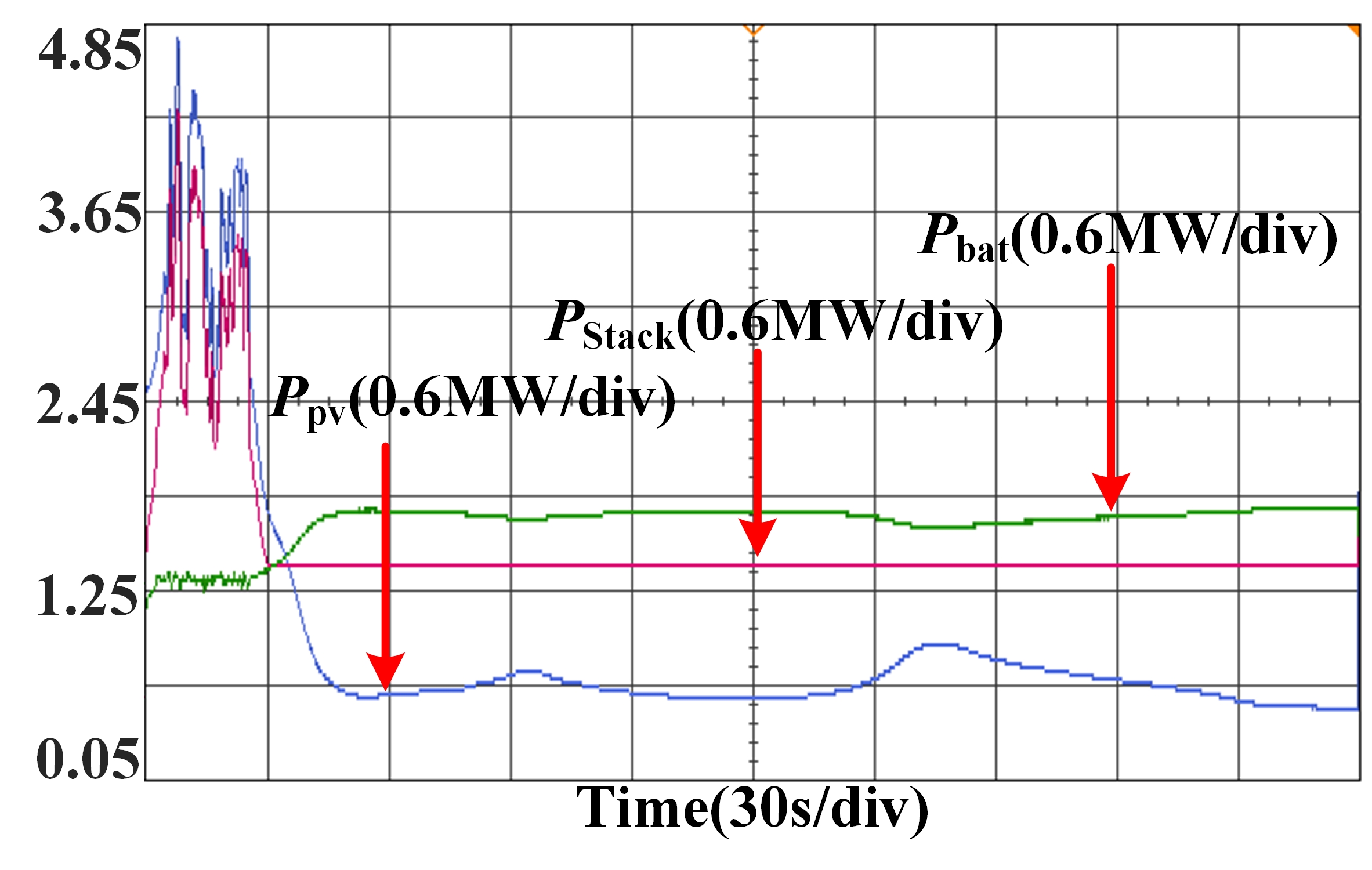}
			%\vspace{-1cm}
			\label{6:e}
		\end{minipage}
	}
	\subfigure[]{
		\begin{minipage}[t]{.3\linewidth}
			\centering
			\includegraphics[scale=0.9]{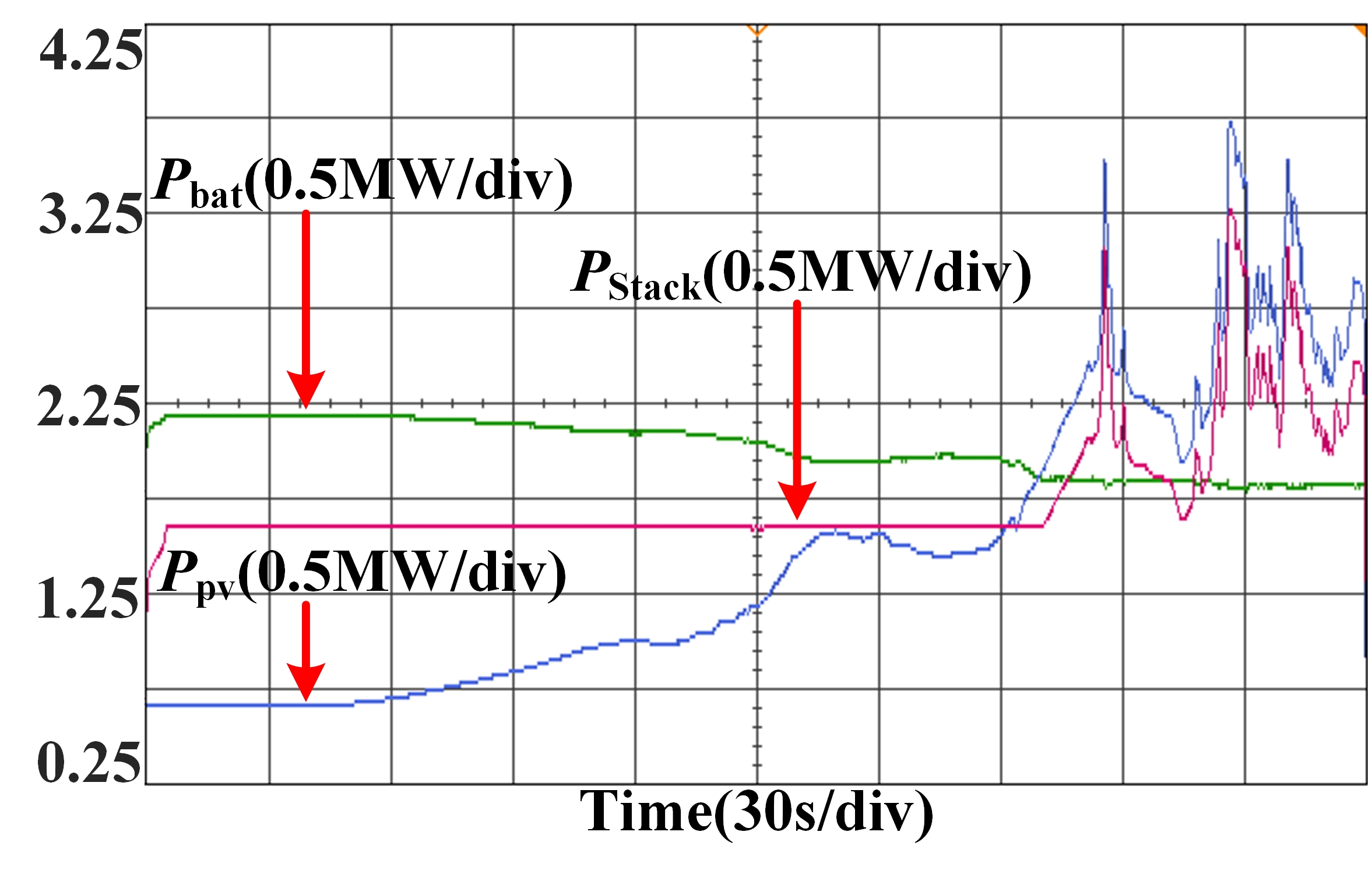}
			%\vspace{-1cm}
			\label{6:f}
		\end{minipage}
	}
	%\vspace{-0.2cm}
	\caption{PV power and electrolytic load in the proposed ESEHE in different scenarios. (a) High-power scenario. (b) Low-power scenario. (c) High-fluctuation scenario. (d) Low-fluctuation scenario. (e) Rapid power-decline scenario. (f) Power-surge scenario.}
	\label{fig:6}%\vspace{-14pt}
\end{figure*}

To validate the megawatt-scale design without prohibitive infrastructure costs, we employ high-fidelity hardware-in-loop (HIL) platform as shown in Fig. \ref{fig:5}. The setup employs a StarSim platform running on an MT8020 real-time simulator with a $1\ \mu$s time step. This high resolution accurately captures the millisecond transient response of the EDL capacitor and the high-frequency switching dynamics of the converter. A TMS320C28346 DSP executes control algorithms at a $10\ \text{kHz}$ sampling rate to ensure the fast response required for stability. Furthermore, an Xilinx XC6SLX16 FPGA generates gate signals for all IGBT switches. Sections \ref{subsec:Grid-Forming Capability Verification} and \ref{subsection:Power Spectrum Analysis of Frequency-Splitting Control} present specific HIL results.

\begin{table}[t]
	\centering\scriptsize
	\renewcommand{\arraystretch}{1.08}  
	\caption{System and Control Parameters in the Simulation} 
	%\vspace{-0.2cm}
	\label{table:2}
	\begin{tabular}{>{\centering\arraybackslash}p{0.05\linewidth}>{\centering\arraybackslash}p{0.04\linewidth}c||>{\centering\arraybackslash}p{0.05\linewidth}>{\centering\arraybackslash}p{0.04\linewidth}c||>{\centering\arraybackslash}p{0.05\linewidth}>{\centering\arraybackslash}p{0.06\linewidth}}  % 使用||创建粗竖线分隔
		\toprule\toprule 
		\multicolumn{1}{c}{\textbf{Stack}} & \textbf{Value} & \textbf{Unit} & 
		\multicolumn{1}{c}{\textbf{Circuit}} & \textbf{Value} & \textbf{Unit} & 
		\multicolumn{1}{c}{\textbf{Control}} & \textbf{Value} \\
		\hline 
		$C_\text{dl}$& 0.02 & F/cm$^\text{2}$ & $V_\text{dc}$& 1000 & V & $K_\text{droop}$& 0.002 \\
		$U_\text{rev}$& 1.228 & V & $C_\text{dc}$& 50 & mF & $K_\text{p}$ & 0.02 \\
		$R_{\text{ohm}}$ & 1.1918 & $\Omega$/cm$^\text{2}$ & $f_{\text{ref}}$ & 50 & Hz & $K_\textit{i}$ & 15 \\
		$I_{\text{exchange}}$ & 0.0015 & A/cm$^\text{2}$ & $L_{\text{buck}}$ & 25 & mH & $K_\text{p1}$ & 0.3 \\
		$N_{\text{cell}}$ & 445 & - & $V_{\text{acref}}$& 35/1 & kV& $K_{\textit{i}\text{1}}$ & 10 \\
		area & 15,000 & cm$^\text{2}$ & $V_{\text{droop}}$ & 0.9697 & V & $K_\text{p2}$ & 0.715 \\
		$K_{\text{act}}$ & 0.1521 & - & $R_{\text{edge}}$ & 0.0009 & $\Omega$ & $K_{\textit{i}\text{2}}$ & 10 \\
		& & & $E_\text{on}$& 40 & mJ & $K_\text{p3}$ & 3 \\
		& & & $E_\text{off}$& 100 & mJ & $K_{\textit{i}\text{3}}$ & 100 \\
		& & & $I_{\text{leak}}$ & 5 & mA & $f_{\text{sw}}$ & 2000\\
		\bottomrule\bottomrule
	\end{tabular}
	%\vspace{-16pt}
\end{table}

\subsection{Grid-Forming Capability Verification}
\label{subsec:Grid-Forming Capability Verification}

\subsubsection{Power-following capability under PV fluctuations}
The directly connected PV-ESEHE system is tested under various PV output scenarios. Results are shown in Fig. \ref{fig:6}.

Figs. \ref{6:a} and \ref{6:b} demonstrate accurate tracking of PV output in both high- and low-power conditions. In high/low fluctuation cases (Figs. \ref{6:c} and \ref{6:d}), the electrolysis stack ramps at $\pm 0.5$ MW/min with zero overshoot and $<50$ ms delay, avoiding high-frequency adjustments and reducing fatigue stress, as discussed in Section \ref{subsec:Frequency-Splitting Power Control Method}. 

Under extreme fluctuations (Figs. \ref{6:e} and \ref{6:f}), for example, such as a 67\% PV drop over 120 s, and the EDL capacitor responds instantaneously through (\ref{14}) to maintain the DC bus voltage $V_\text{dc}$. The ES switches to discharge mode according to (\ref{8}), (\ref{9}) to compensate for the power deficit and maintain stable electrolysis output. Finally, the electrolytic load is down-adjusted within a rate range of $\pm 0.5$ MW/min to alleviate electrode degradation.

These results confirm that the proposed ESEHE autonomously stabilizes the off-grid ReP2H system, supports GFL PV operation, and sustains continuous hydrogen production even under high volatility.

\begin{figure*}[htbp]
	\centering
	\subfigure[]{
		\begin{minipage}[t]{.3\linewidth}
			\centering
			\includegraphics[scale=0.94]{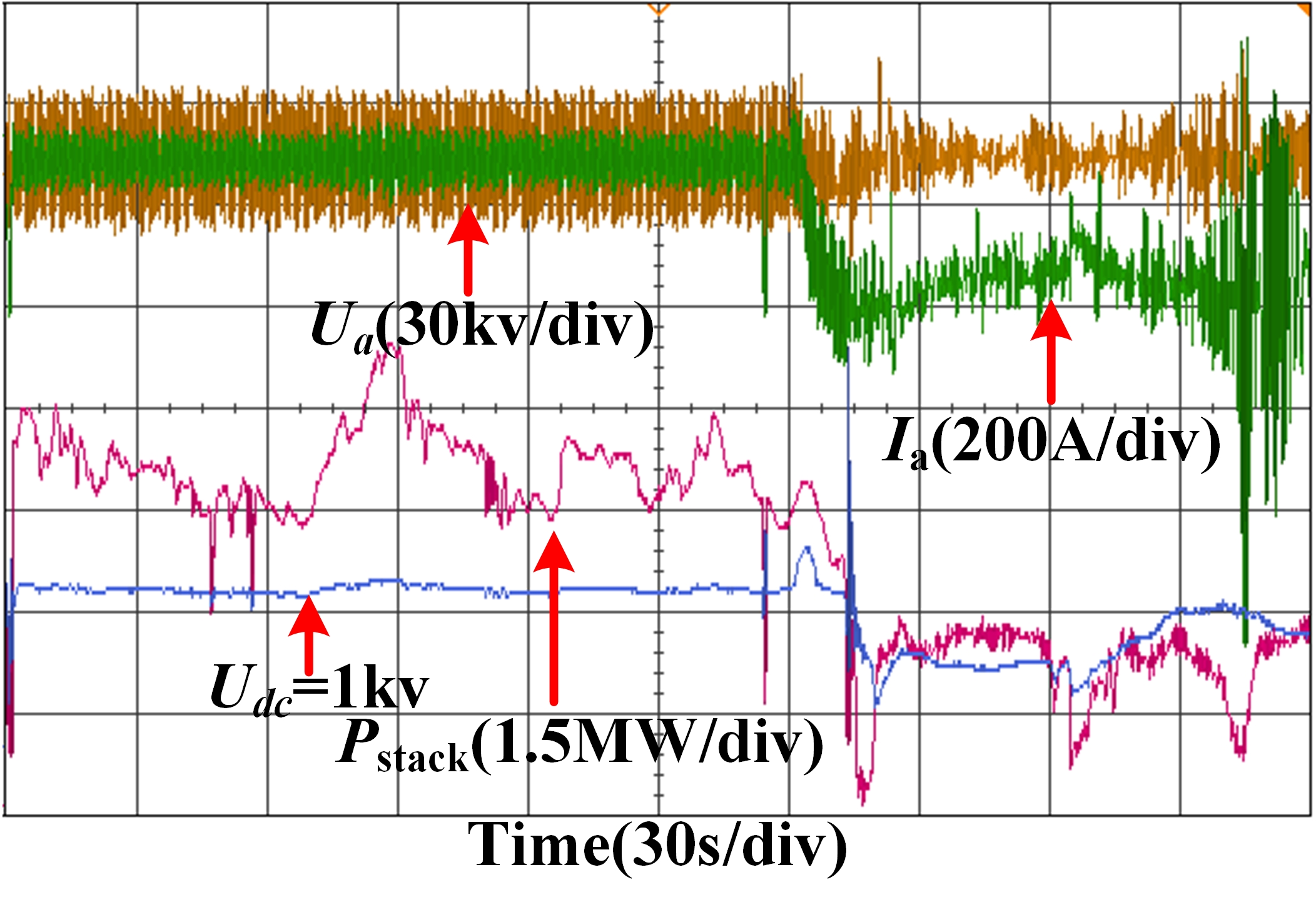}
			\label{7:a}
		\end{minipage}
	}
	\subfigure[]{
		\begin{minipage}[t]{.3\linewidth} 
			\centering
			\includegraphics[scale=0.94]{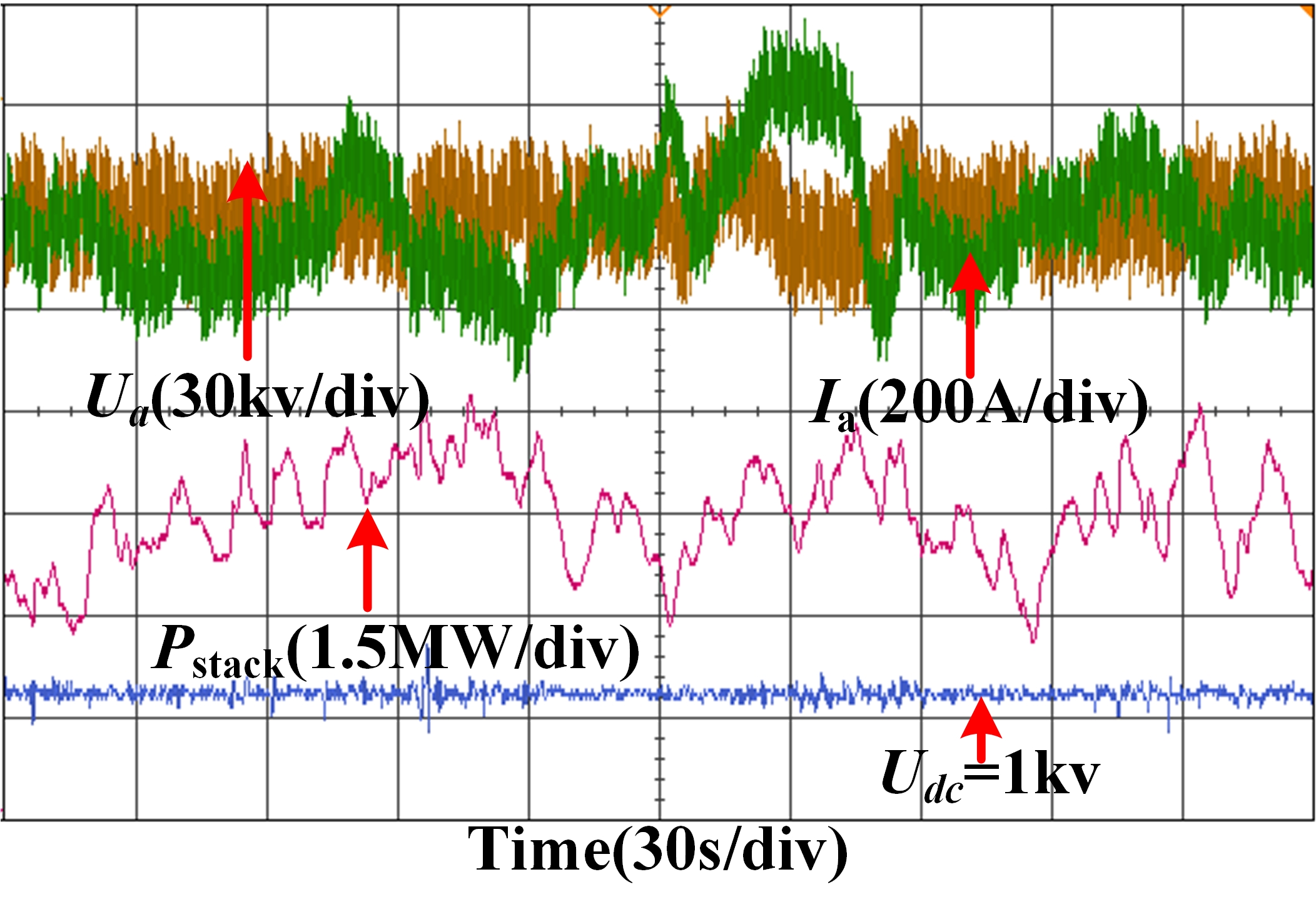}
			\label{7:b}
		\end{minipage}
	}
	\subfigure[]{
		\begin{minipage}[t]{.3\linewidth}
			\centering
			\includegraphics[scale=0.94]{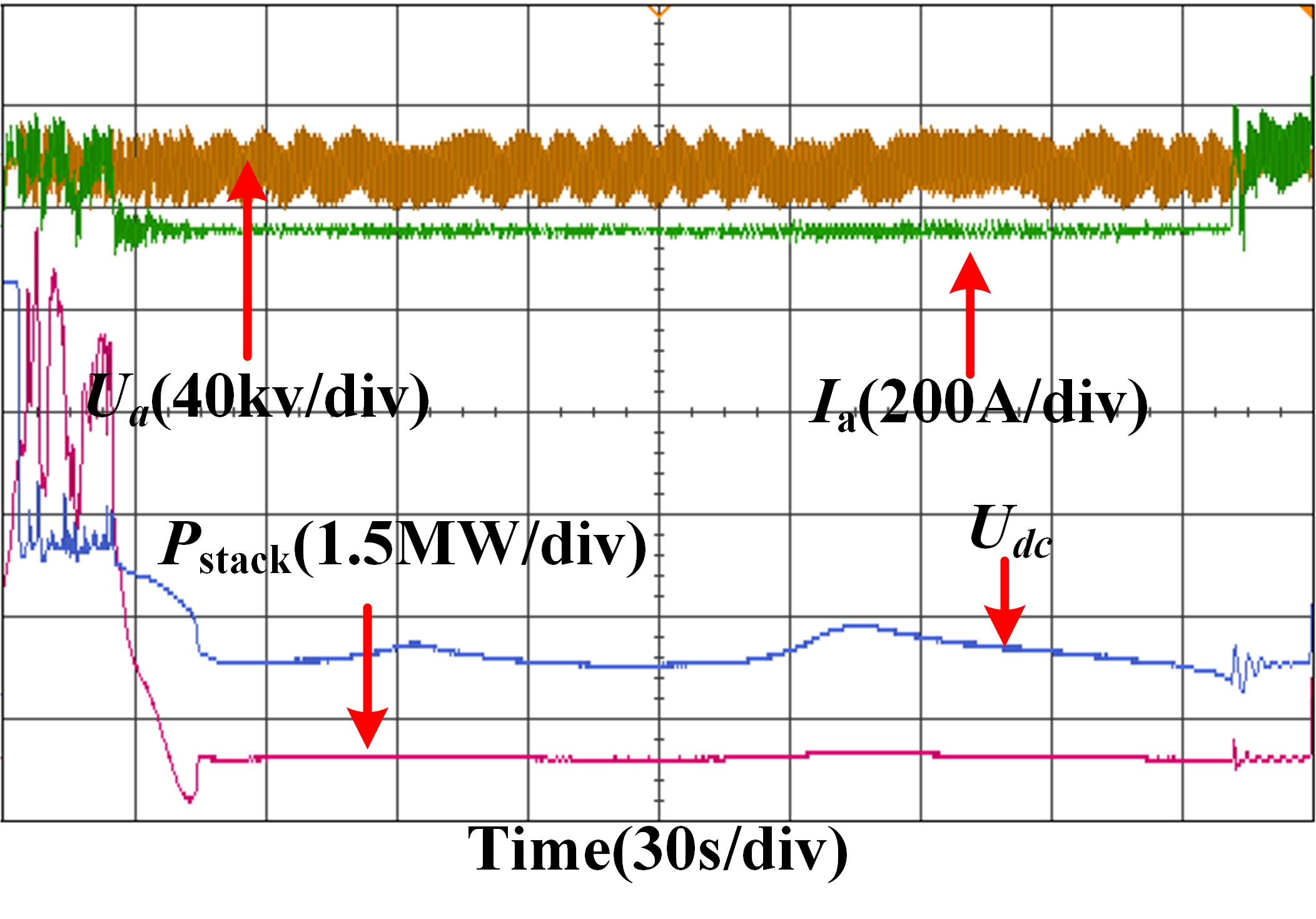}
			\label{7:c}
		\end{minipage}
	}
	\subfigure[]{
		\begin{minipage}[t]{.3\linewidth} 
			\centering
			\includegraphics[scale=0.94]{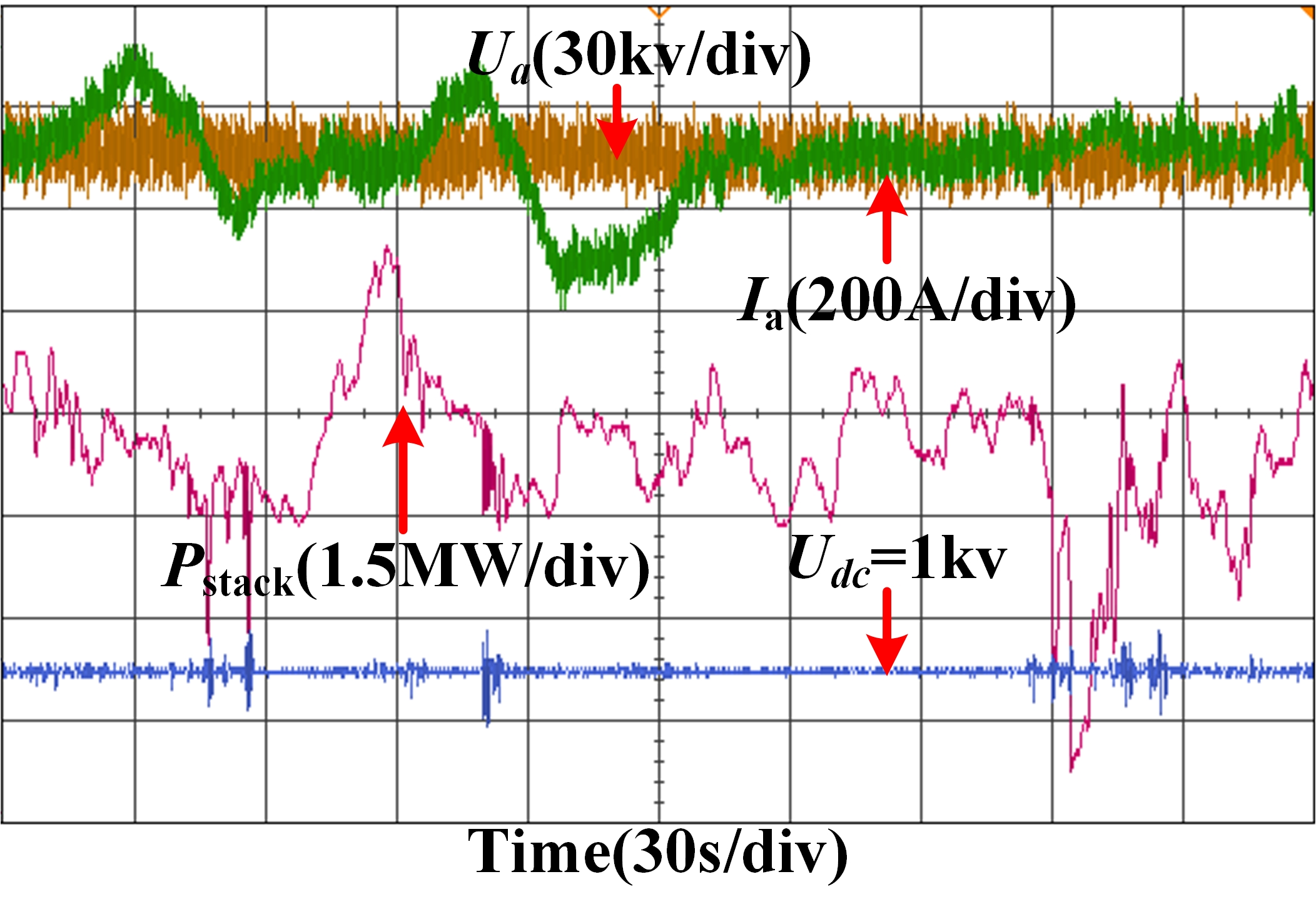}
			\label{7:d}
		\end{minipage}
	}
	\subfigure[]{
		\begin{minipage}[t]{.3\linewidth} 
			\centering
			\includegraphics[scale=0.94]{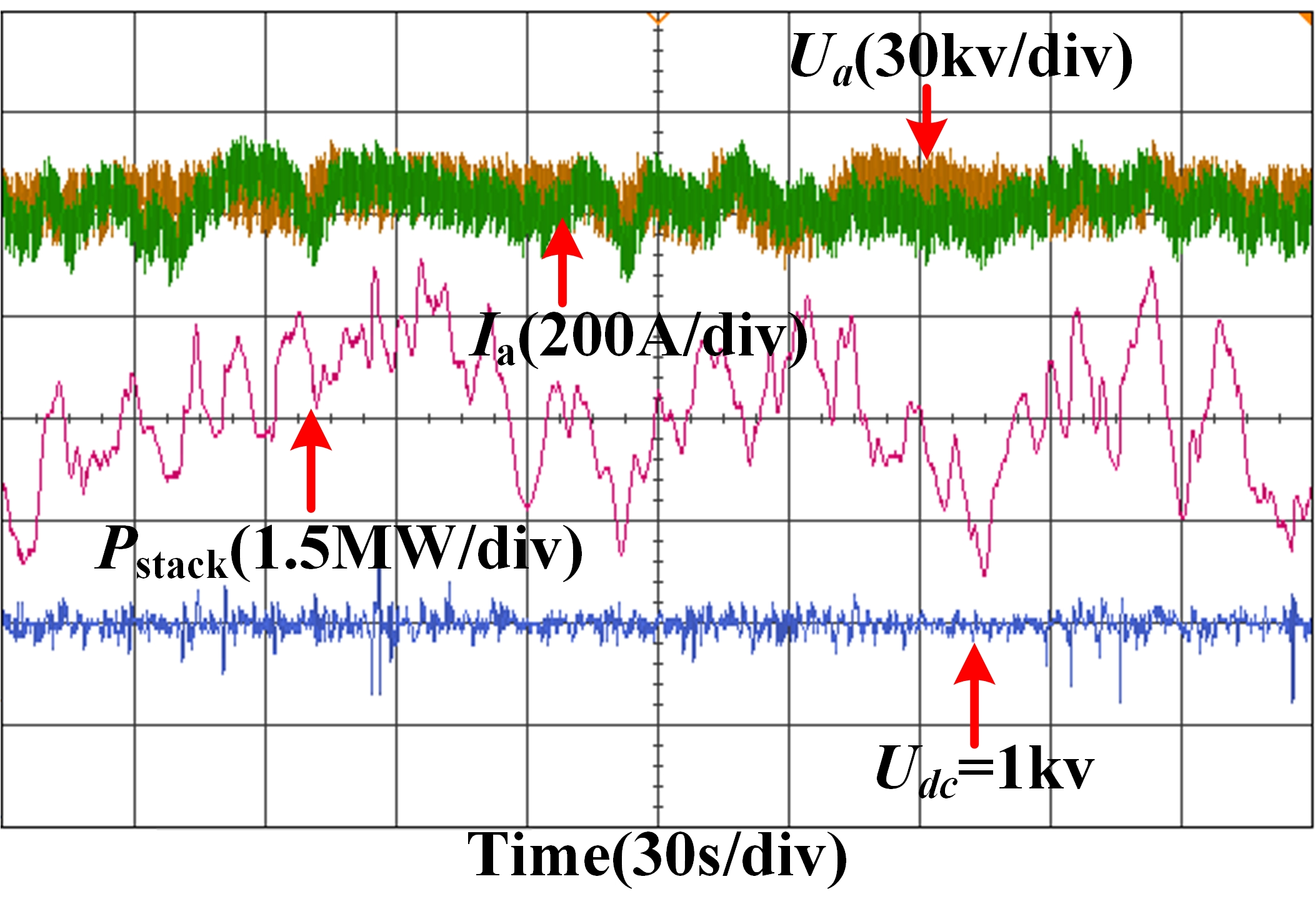}
			\label{7:e}
		\end{minipage}
	}
	\subfigure[]{
		\begin{minipage}[t]{.3\linewidth} 
			\centering
			\includegraphics[scale=0.94]{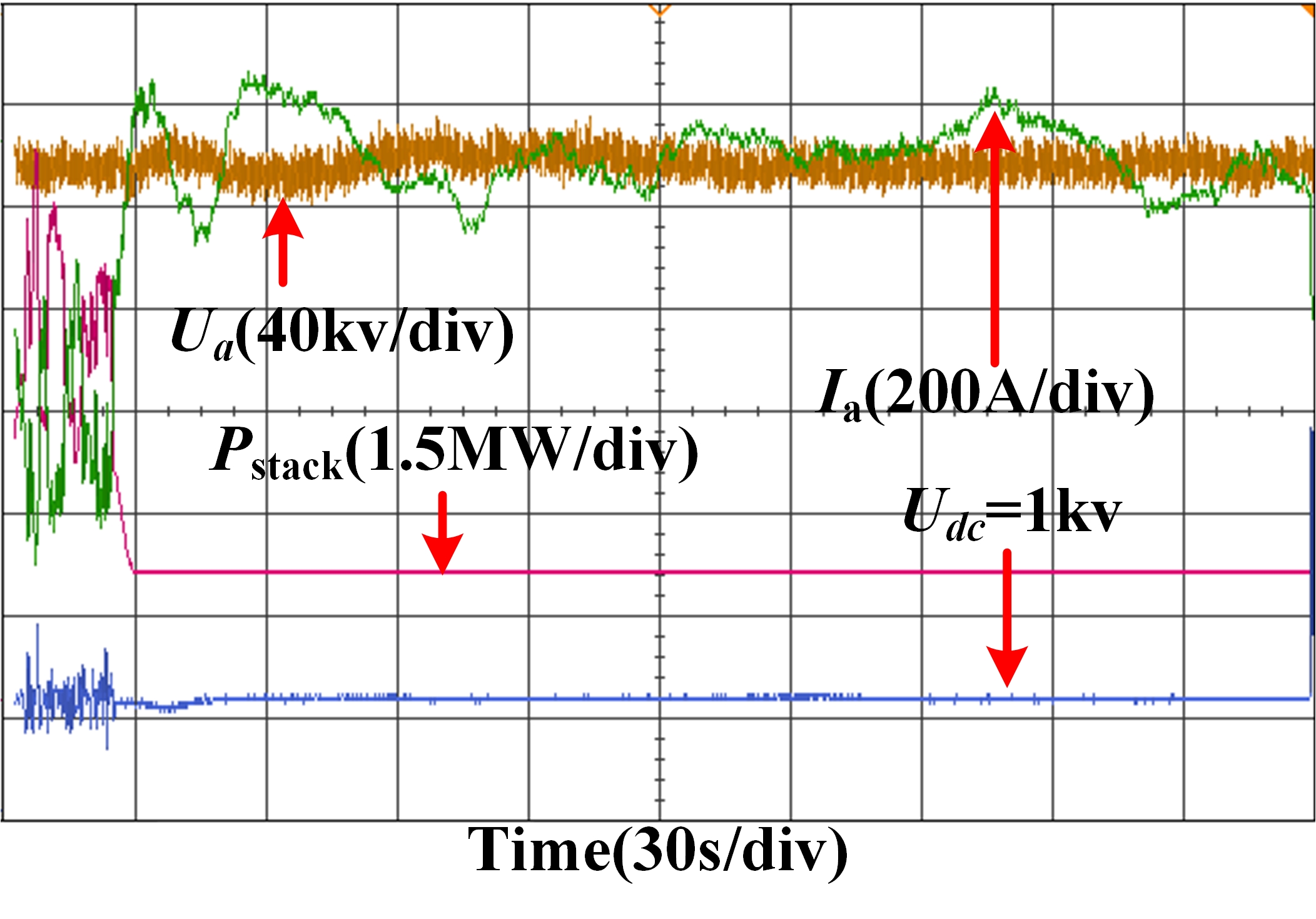}
			\label{7:f}
		\end{minipage}
	}
	%\vspace{-0.2cm} 
	\caption{Comparative responses of AC bus voltage/current, DC bus voltage, and electrolysis power between conventional GFM HE and the proposed ESEHE under different power levels and fluctuation scenarios. (a) Conventional GFM HE without ES under high-volatility scenario. (b) Conventional GFM HE without ES under high-power scenario. (c) Conventional GFM HE without ES under rapid power-decline scenario. (d) The proposed ESEHE under high-volatility scenario. (e) The proposed ESEHE under high-power scenario. (f) The proposed ESEHE under rapid power-decline scenario.}
	\label{fig:7}
\end{figure*}

\begin{figure*}[htbp]
	\centering
	\subfigure[]{
		\begin{minipage}[t]{.3\linewidth} 
			\centering
			\includegraphics[scale=0.94]{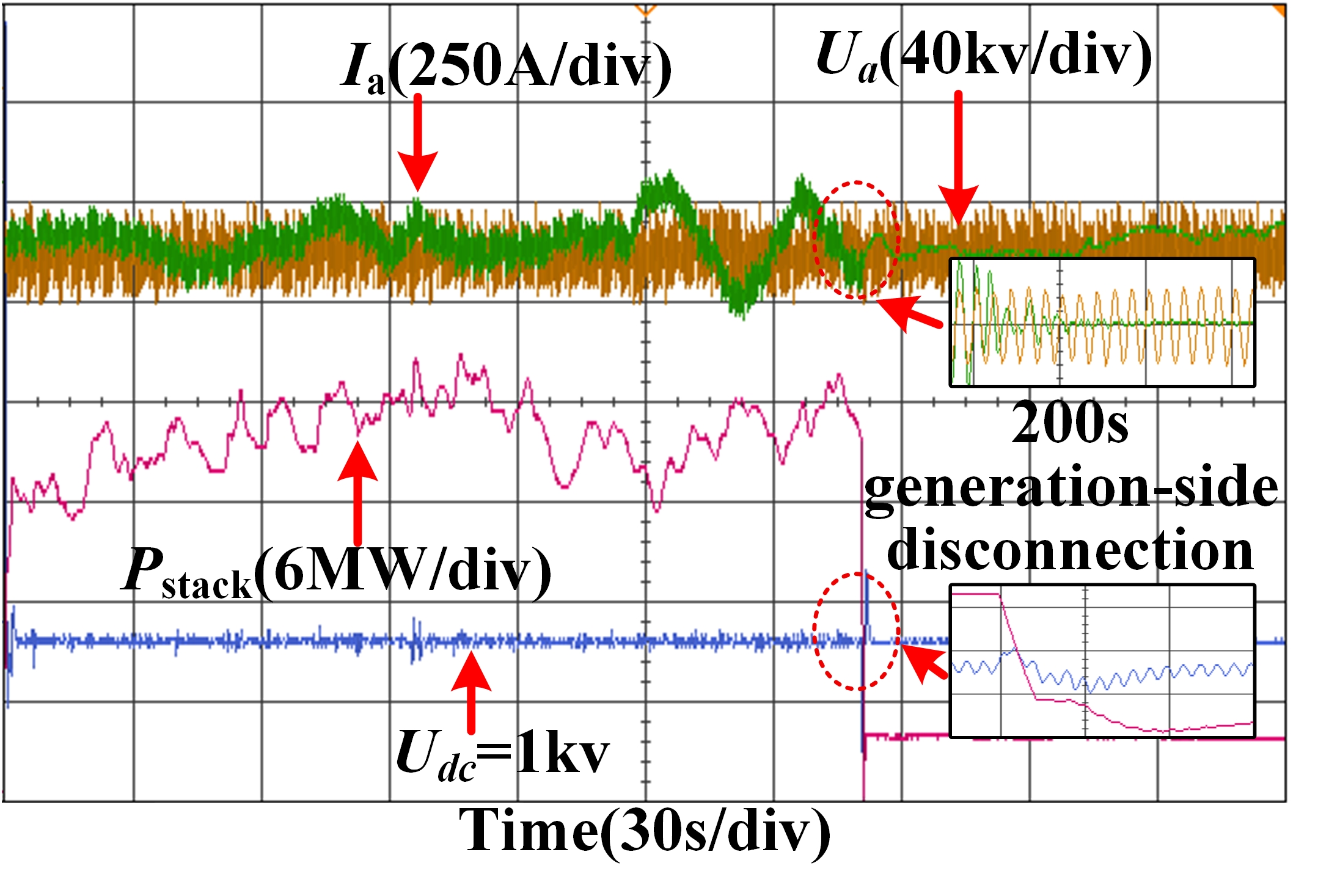}
			\label{8:a}
		\end{minipage}
	}
	\subfigure[]{
		\begin{minipage}[t]{.3\linewidth}
			\centering
			\includegraphics[scale=0.94]{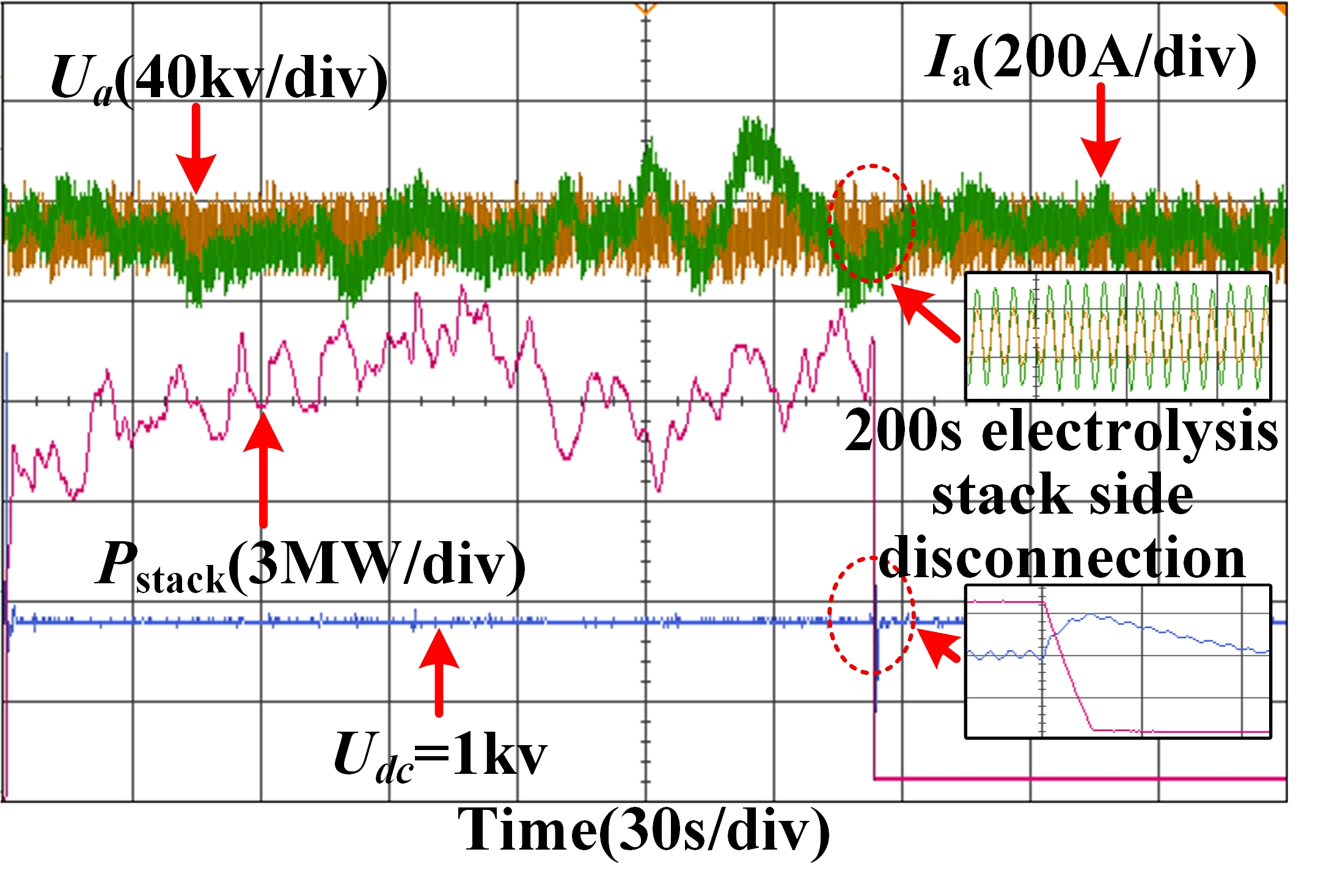}
			\label{8:b}
		\end{minipage}
	}
	\subfigure[]{
		\begin{minipage}[t]{.3\linewidth} 
			\centering
			\includegraphics[scale=0.94]{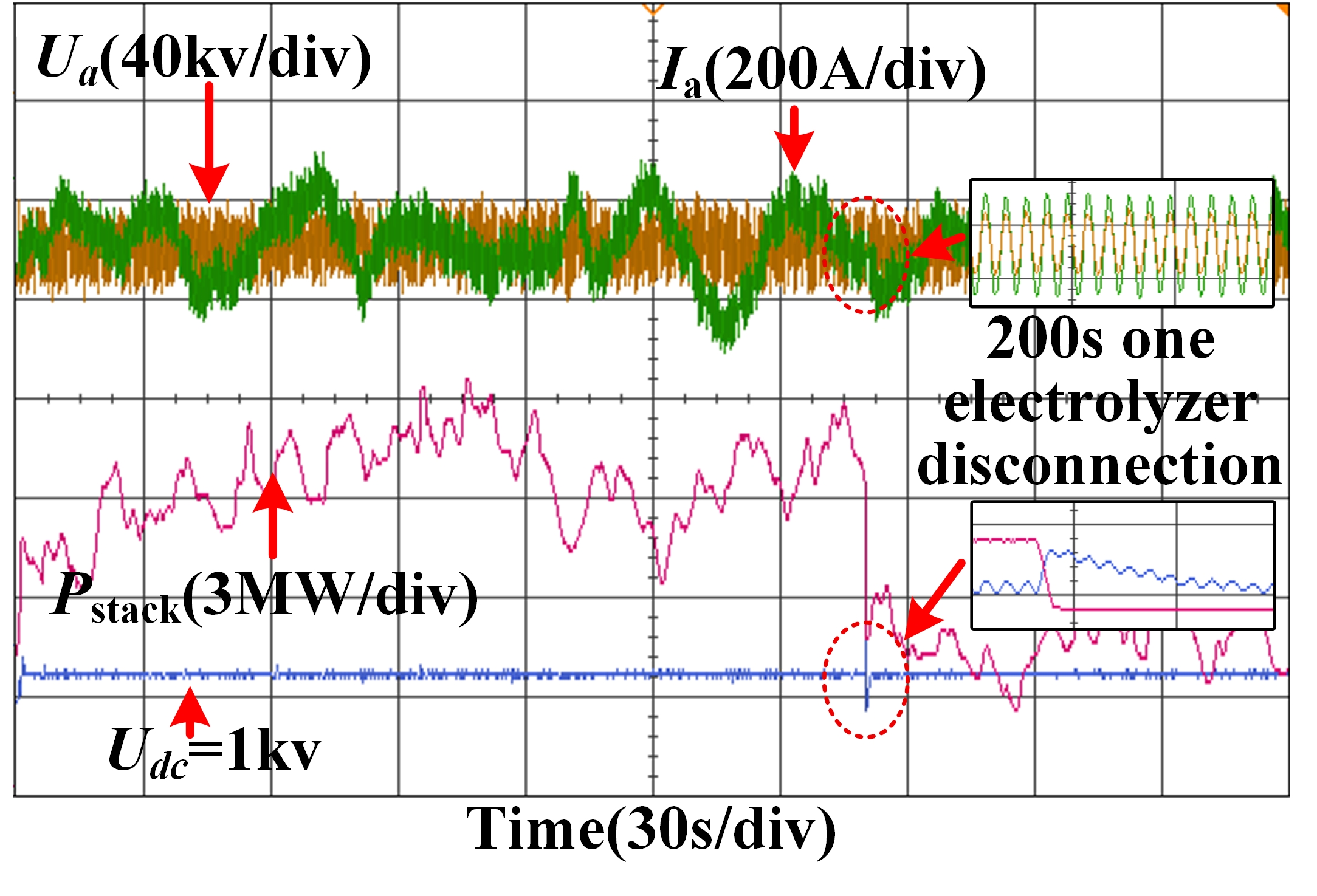}
			\label{8:c}
		\end{minipage}
	}
	
	\caption{Responses of AC bus voltage/current, DC bus voltage, and electrolysis power for the proposed ESEHE under $N$-1 contingency. (a) Generation-side disconnection. (b) Electrolysis stack disconnection. (c) Two ESEHEs operate in parallel when one electrolyzer is disconnected.}
	\label{fig:8}
	%\vspace{-0.5cm} 
\end{figure*}

\begin{table}[t]
	\centering
	\renewcommand{\arraystretch}{1.2}
	\caption{Performance Comparison Between Conventional GFM Hydrogen Electrolyzer without ES and the Proposed ESEHE}
	%\vspace{-0.3cm} 
	\label{table:3} 
	\begin{tabular}{@{}>{\centering\arraybackslash}m{0.35\columnwidth}>{\centering\arraybackslash}m{0.32\columnwidth}>{\centering\arraybackslash}m{0.2\columnwidth}@{}}
		\toprule 
		\toprule
		\textbf{Performance Indicators} & 
		\textbf{Conventional GFM Electrolyzer (w/o ES)}& 
		\textbf{Proposed ESEHE}\\
		\midrule 
		AC voltage range (p.u.)& $\pm$9.5\%-$\pm$15\%& $\pm$3\%-$\pm$5\%\\frequency deviation (Hz)& $\pm$1.7-$\pm$2.1& $\pm$0.5-$\pm$0.8\\
		Power fluctuation response & Poor & Good \\
		\bottomrule 
		\bottomrule 
	\end{tabular}
	%\vspace{-14pt}
\end{table}

\subsubsection{Comparative analysis of grid-forming HEs with and without ES}

Fig. \ref{fig:7} compares the proposed ESEHE with the conventional GFM HE proposed in \cite{Grid-FormingServicesFromHydrogenElectrolyzers}.

For the conventional HE without ES (Figs. \ref{7:a}--\ref{7:c}), high fluctuation, high-power, and sharp power-drop cases (e.g., 67\% drop) cause bus voltage instability and power deviation due to insufficient energy buffering, even with GFM control. 

In contrast, the proposed ESEHE (Figs. \ref{7:d}--\ref{7:f}) maintains bus voltage stability ($\pm 2\%$ nominal value) and precisely tracks PV fluctuations ($\pm 0.5$ MW/min).  %through frequency-splitting control and adaptive grid-forming control. 
Additionally, by directly connecting the ES to the DC bus, the proposed ESEHE reduces energy flow across the AC/DC conversion links, shortens transient response time, and demonstrates both topological and control advantages. Quantitative metrics are given in Table \ref{table:3}.

\subsubsection{$N$-1 contingency verification}
\label{sec:n-1}

Figs. \ref{8:a}--\ref{8:c} show the system’s transient performance under three $N$-1 faults: loss of generation, loss of electrolysis stack, and loss of an electrolyzer in a multi-unit system. Detailed analysis is as follows: 

a. In Fig. \ref{8:a}, generation is disconnected at 200 s.  The resulting power deficit causes the DC bus voltage $V_\text{dc}$ to drop. In accordance with (\ref{3}), the VSM controller raises the reference angular frequency $\omega^{*}$, driving the ESEHE to switch instantly to \emph{mode h} in Table \ref{table:1} and minimize the electrolytic load. Concurrently, the ES discharges through its DC/DC converter to maintain the AC and DC voltages within $\pm 3\%$ and $\pm 2\%$ of nominal values, respectively. The system returns to steady state within 150 ms, highlighting rapid transient recovery.

b. In Fig. \ref{8:b}, as the stack is disconnected at 200 s, the DC bus voltage $V_\text{dc}$ rises. In response,  $\textit{P}_\text{ref}$ in (\ref{4}) instantaneously drops to 0 to alleviate the voltage rise in the DC bus, PV power is redirected to the ES, and ESEHE switches to \emph{mode e} in Table \ref{table:1}. Moreover, the EDL capacitor instantly damps the ensuing voltage spikes via (\ref{14}) to ensure that the AC/DC voltage is maintained between $\pm 5\%$ and $\pm 3\%$ of nominal values. Stabilization is achieved 150 ms after the fault.

c. Fig. \ref{8:c} shows one of two parallel ESEHEs, powered by two GFL PV units, being disconnected at 200 s. After the fault, AC and DC voltages briefly rise (within 2\% and 3\% of nominal values, respectively) due to surplus energy, but are stabilized quickly. During this process, the EDL capacitor instantly absorbs high-frequency components and buffers the transient voltage changes. The ES charges to absorb excess power and respond to mid-frequency power fluctuations. The remaining ESEHE adjusts the $\textit{P}_\text{ref}$ according to (\ref{4}) to continue following PV output, thus maintaining overall stability.

\subsection{Power Spectrum Analysis of Frequency-Splitting Control}
\label{subsection:Power Spectrum Analysis of Frequency-Splitting Control}

\begin{figure}[!ht]
	\centering
	\includegraphics[width=0.70\linewidth]{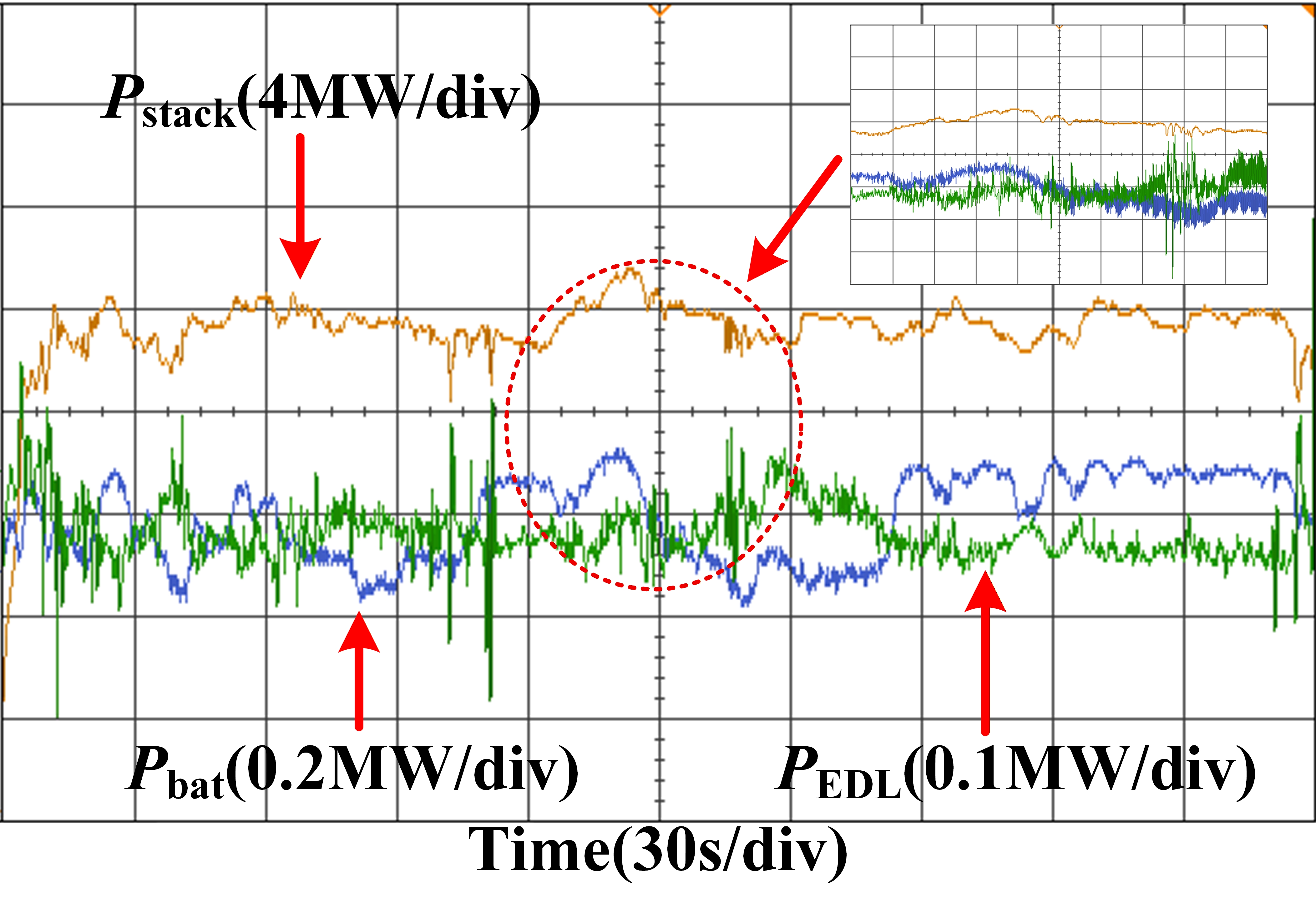}
	%\vspace{-6pt}
	\caption{Power responses of EDL capacitor, ES, and electrolytic load under  high-volatility scenario.}
	\label{fig:d9}
\end{figure}

Fig. \ref{fig:d9} shows the transient power responses of system components under high volatility conditions to verify the frequency splitting control strategy in the time domain. The proposed controller decomposes the variable renewable power input, allowing the ES branch to absorb mid-frequency fluctuations while the EDL capacitor manages instantaneous high-frequency spikes. As a result, the electrolytic power remains smooth and follows the low-frequency command. This behavior confirms that the strategy isolates the stack from dynamic stress.

\begin{figure}[!ht]
	\centering
	\includegraphics[width=0.78\linewidth]{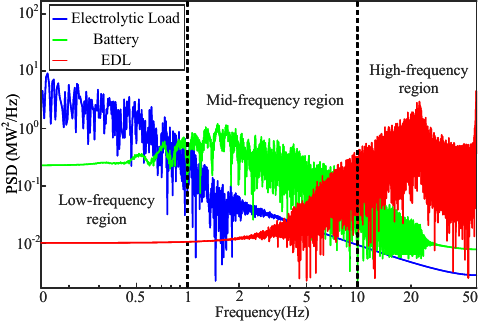}
	%\vspace{-6pt}
	\caption{Power spectrum of the electrolytic load, EDL capacitor, and ES for the proposed ESEHE under frequency-splitting control.}
	\label{fig:9}
\end{figure}

Fig. \ref{fig:9} presents the effect of the frequency-splitting control described in Section \ref{subsec:Frequency-Splitting Power Control Method} under a fluctuating PV scenario (same as Fig. \ref{6:c}). Power spectral density (PSD) analysis shows that the EDL capacitor leverages its millisecond-scale response to absorb high-frequency components, limiting the rate of change of electrode current to 0.1 A/$\mu$s, thus suppressing fluctuations in electrode overvoltage. The ES mitigates stress on electrodes and stack structures by absorbing mid-frequency components via band-pass filtering. Meanwhile, the electrolytic load seamlessly tracks low-frequency power with a constrained ramp rate ($\pm 0.5 $ MW/min, in Fig. \ref{6:c}), thus avoiding accelerating fatigue cracks and electrode degradation of the stack.

Moreover, by confining the ES power to the intermediate-frequency band, the strategy avoids medium- and high-frequency charge/discharges, % that cause nonuniform lithium-ion deposition, 
thereby reducing %cycling-induced
degradation and extending battery life. The corresponding economic benefit is assessed in Section \ref{ES degradation loss}.

\subsection{Verification of SOC Equalization for Multiple ESEHEs }
\label{subsec:3}

\begin{figure}[!ht]
	\centering
	\includegraphics[width=0.70\linewidth]{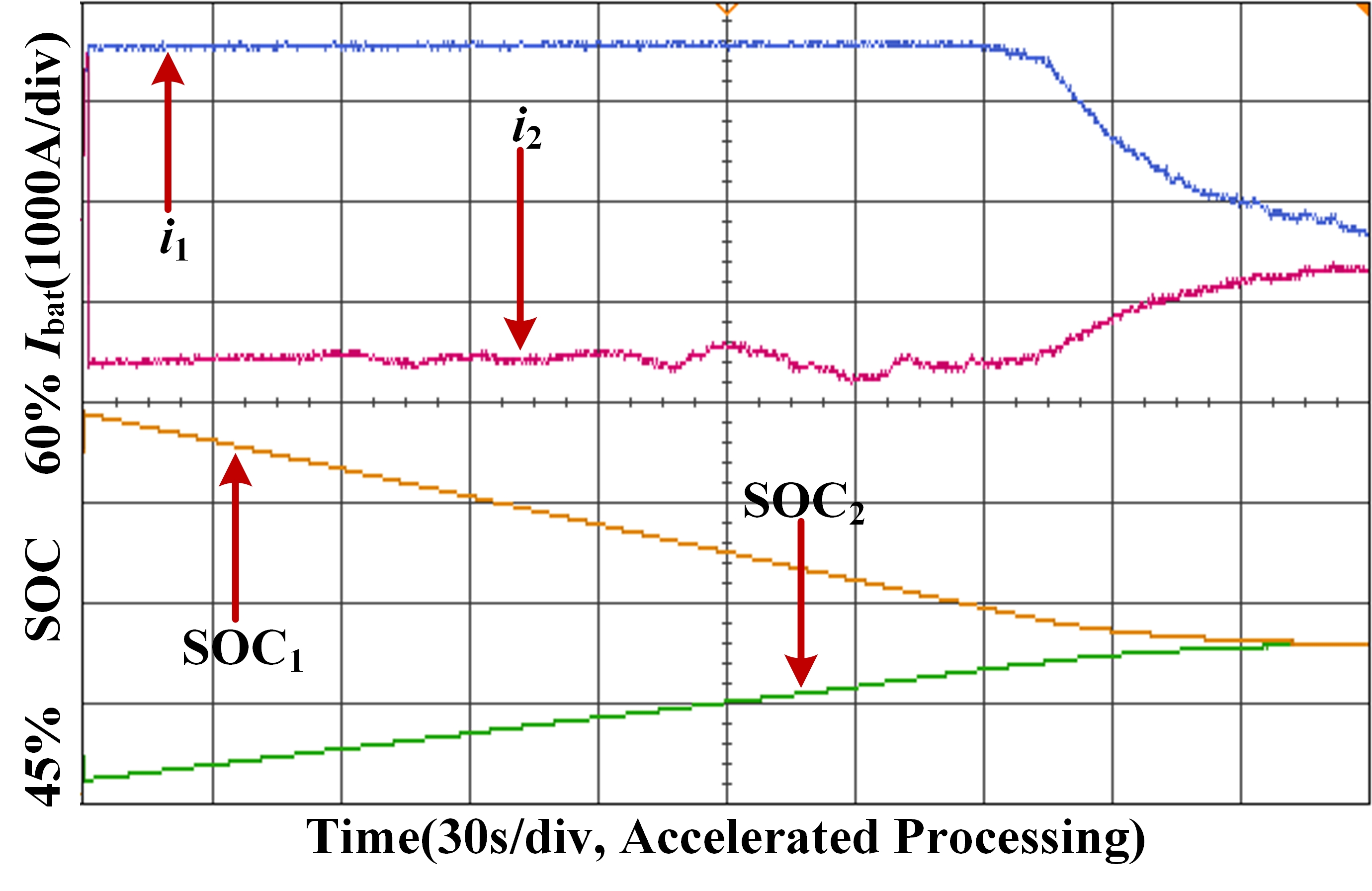}
	%\vspace{-6pt}
	\caption{SOC of the ES branches of two ESEHEs under PV fluctuations.}
	\label{fig:10}
	%\vspace{-16pt}
\end{figure}

Following the setup in Section \ref{sec:n-1} and Fig. \ref{fig:4}, two ESEHEs are connected to two PV units. The initial SOCs of the ES batteries are 55\% and 45\%. The upper and lower limits of SOC are set to 90\% and 10\%. The PV output follows the same profile as in Fig. \ref{6:a}. A 3-hour PV dataset is time-compressed and simulated in 5 minutes, and the simulation gives the SOC trajectories as shown in Fig. \ref{fig:10}.

The SOCs converge within 0.5\% in 300 s of simulation. Given the real-time compression ratio of 36:1, this corresponds to an equalization rate of 1.83\%/min in real time. During PV fluctuations, the discharging power shares of the higher- and lower-SOC ESEHEs are 63\% and 71\%, respectively, driven by the dynamic droop control in Section \ref{subsec:Equalization Control Strategy}. After stabilization, the SOC differences between the two ESEHEs are below 1\%, preventing deep discharging of lower-SOC batteries.
Without a centralized controller, this strategy can also support scalable deployment of large off-grid ReP2H systems.

\subsection{Techno-Economic Benefit Analysis}
\label{subsec:Techno-Economic Benefit Analysis}

This section compares the techno-economic performance between the proposed ESEHEs and a traditional off-grid ReP2H system based on conventional electrolyzers and a centralized GFM ES plant. Energy losses and ES degradation are quantified. The rating of the electrolysis stack is set to 5 MW, and the ES branch has a capacity of 0.8 MWh.

\subsubsection{Converter losses reduction}

Converter losses arise from the AC/DC interface and DC/DC converters for both the electrolysis stack and ES branches, including conduction and switching losses.
From \cite{PowerLossandEfficiencyCalculation}, we model the converter losses:
\begin{gather}
	P_{\text{cond}} = V_{\text{ce}} \cdot I_{\text{avg}} + R_{\text{ce}} \cdot I_{\text{rms}}^2,
	\label{42} \\
	P_{\text{sw}} = \left( E_{\text{on}} + E_{\text{off}} \right) \cdot f_{\text{sw}},
	\label{43}
\end{gather}
where $\textit{V}_\text{ce}$ is the on-state voltage drop; $\textit{I}_\text{avg}$ is the average current per module; $\textit{R}_\text{ce}$ is the conduction marginal resistance; $\textit{I}_\text{rms}$ is the RMS current (equal to $\textit{I}_\text{avg}$ under steady hydrogen production; $\textit{E}_\text{on}$ and $\textit{E}\text{off}$ are the turn-on and turn-off switching energy per cycle; and $\textit{f}_\text{sw}$ is the switching frequency. The selection of converter parameters is based on the converter loss models from datasheets \cite{PowerLossandEfficiencyCalculation, LIReal-time}.

In the ESEHE, the ES branch discharges to support hydrogen production when renewable power supply drops, as in Fig. \ref{6:e}. The energy passes from the ES battery to the stack through only two DC/DC converters. In contrast, the conventional system uses separate converters for the HE and ES plant, requiring two AC/DC and two DC/DC conversions, thus causing higher losses. On average, each AC/DC converter contributes 5.82 kW conduction loss and 280 W switching loss, totaling 11.64 kW. By eliminating these, the ESEHE achieves an average efficiency gain of 0.233\%.

\subsubsection{ES degradation reduction}
\label{ES degradation loss}

\begin{figure}[!t]
	\centering
	\includegraphics[width=0.75\linewidth]{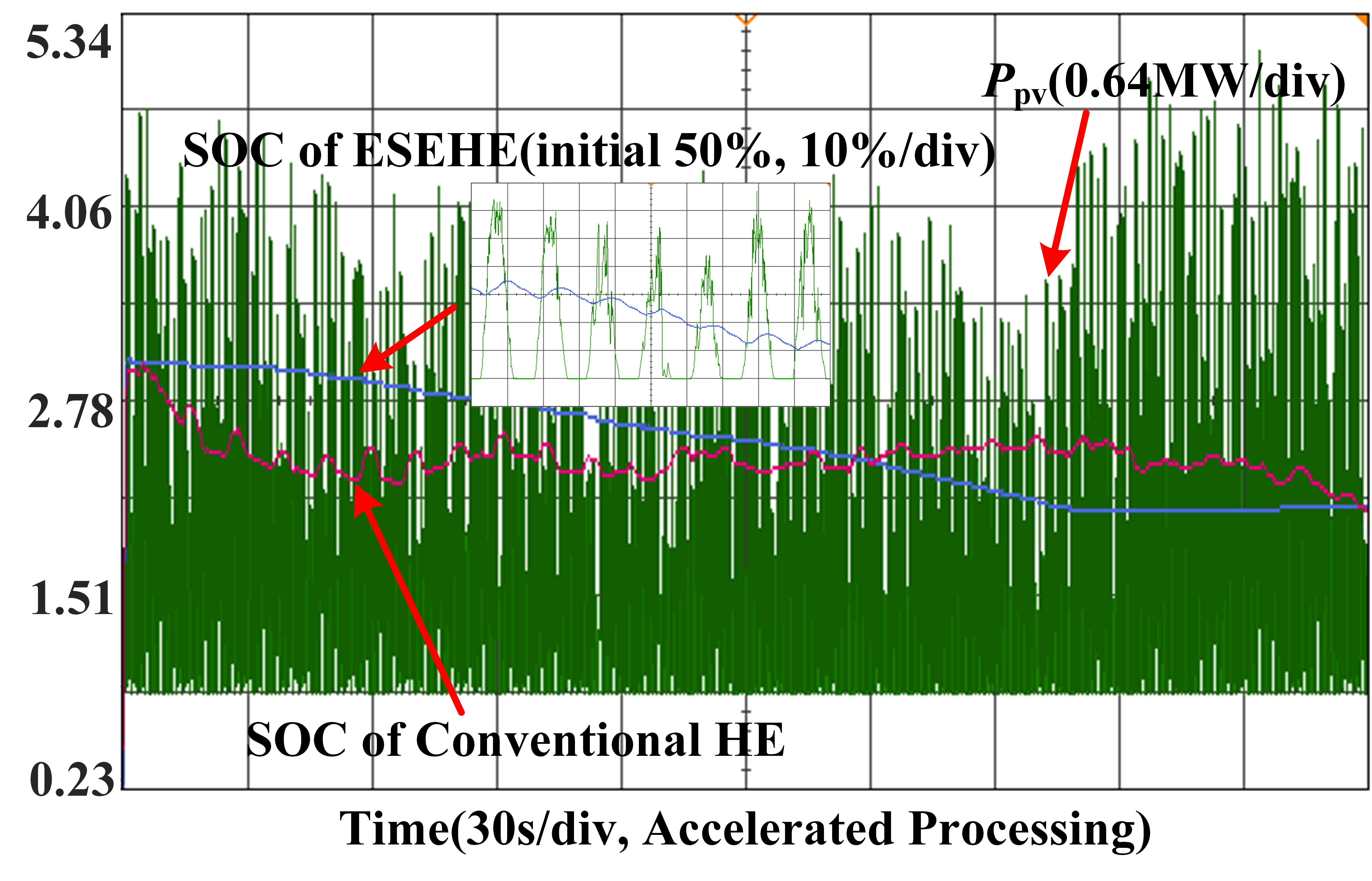}
	%\vspace{-0.4cm}
	\caption{Annual PV power and SOC variation in techno-economic analysis.}
	\label{fig:11}
	%\vspace{-0.2cm}
\end{figure}

ES degradation is estimated from the complete charge-discharge cycle count using rain-flow decomposition of SOC data {($\textit{$\gamma$}\text{SOC}_\text{k}$, $\textit{N}_\text{k}$)}, where:
\begin{equation}
	\gamma \text{SOC}_\text{k} = \text{SOC}_\text{max,k} - \text{SOC}_\text{min,k},
	\label{44}
\end{equation}
where $\textit{$\gamma$}\text{SOC}_\text{k}$ denotes the depth of the $\textit{k}$-th cycle; and $\textit{N}_\text{k}$ is the equivalent cycle count. Then, we apply Miner's linear damage accumulation theory to assess the battery cycle life:  
\begin{equation}
	D_{\text{total}} = \sum_{\text{k}=1}^{n} \frac{N_\text{k}}{N_f(\gamma \text{SOC}_\text{k})},
	\label{45}
\end{equation}
where $\textit{N}_{\text{f}}(\textit{$\gamma$}\text{SOC}_\text{k}) = \textit{$\alpha$}(\textit{$\gamma$}\text{SOC}_\text{k})^{-\textit{{$\beta$}}}$ is the cycle-to-failure at depth $\textit{$\gamma$}\text{SOC}_\text{k}$; $\alpha$ and $\beta$ are constants. As per Fig. \ref{fig:11}, since the ES operates under shallow cycling conditions, $\alpha$ and $\beta$ are set to $4,000$ and $0.47$ based on battery empirical data \cite{zhu2024exploring, PowerLossandEfficiencyCalculation}. Failure occurs when $\textit{D}_\text{total}$ $\geq$ 1. 

Based on 2022 PV data from Inner Mongolia (Fig. \ref{fig:11}) and the EMS framework in \cite{zhu2024exploring}, simulations show that the ESEHE exhibits an average daily SOC swing of 5\% occurring twice per day, equivalent to approximately 730 full cycles per year. According to (\ref{40}), this corresponds to a cycle life $\textit{N}_{\text{f}}$ of 16,351 cycles at a 5\% SOC depth. As a result, the annual ES degradation rate is 4.5\%, lower than the 5.6\% observed in the conventional HE system (purple curve in Fig. \ref{fig:11}) \cite{PowerLossandEfficiencyCalculation}, extending battery lifetime by 19.7\%.

When the ES capacity decreases to 80\%, it is considered end-of-life. Assuming a 20-year lifespan, the ES in a conventional ReP2H system has a 5.6\% annual degradation rate and needs battery replacement every 3.57 years, whereas the proposed ESEHE extends this interval to 4.44 years, postponing the first replacement by 0.87 years. Over the lifecycle, the conventional ES needs five replacements, while the ESEHE needs only four, resulting in a cumulative delay of about 2.2 years. Considering cost parameters for a 0.8 MWh system (replacement 900 CNY/kWh; recycling revenue 150 CNY/kWh \cite{zhu2024exploring}), the ESEHE achieves total lifecycle savings of 0.6 million CNY.

To validate economic robustness against model simplifications and environmental uncertainties, we conduct a quantitative sensitivity analysis. We examine four parameters within a 20\% variation range, including battery cycle life, operating temperature, switching energy losses, and SOC estimation error.
	
\begin{figure}[!ht]
	\centering
	\includegraphics[width=0.8\linewidth]{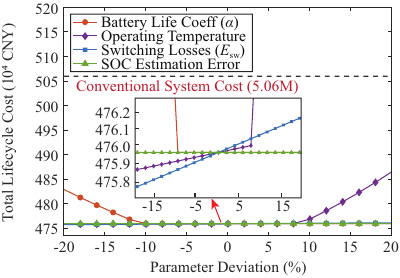}
	%\vspace{-0.4cm}
	\caption{Sensitivity analysis of the total life-cycle cost (TLC) to variations in key parameters.}
	\label{fig:14}
	%\vspace{-0.2cm}
\end{figure}

Fig. \ref{fig:14} shows the impact of these variables on the total life-cycle cost. The system proves most sensitive to battery cycle life and temperature. A decline in battery quality or a sustained rise in temperature causes a non-linear cost increase by accelerating aging and necessitating frequent replacements. In contrast, the flat curve for switching losses indicates that minor efficiency drifts have a negligible impact on overall costs. Similarly, the unnoticeable variation caused by the SOC estimation error confirms that the proposed equalization strategy is robust, as deviations do not accelerate degradation.

Crucially, the total life-cycle cost of the proposed ESEHE remains consistently below the conventional system benchmark of 5.06 million CNY. This confirms the solution maintains economic superiority despite potential parameter variations.

\begin{table}[t]
	\centering
	\renewcommand{\arraystretch}{1.2} 
	\caption{Comparison of Converter Costs Between Conventional Hydrogen Electrolyzer and the Proposed ESEHE}
	%\vspace{-0.2cm}
	\label{table:5}
	\begin{tabular}{@{}>{\centering\arraybackslash}m{0.3\columnwidth}>{\centering\arraybackslash}m{0.28\columnwidth}>{\centering\arraybackslash}m{0.28\columnwidth}@{}}
		\toprule 
		\toprule
		\textbf{Block} & \textbf{Conventional HE} & \textbf{Proposed ESEHE} \\
		\midrule 
		Electrolyzer AC/DC		&  5.15 MVA	&  5.22 MVA\\
		ES AC/DC 				&  0.82 MVA	& -- \\ 
		Stack DC/DC  			&  5.15 MVA	&  5.15 MVA\\
		Battery DC/DC			&  0.82 MVA	&  0.82 MVA\\ \midrule 
		Total AC/DC capacity             
		& 	5.97 MVA  
		&  	5.22 MVA    \\
		Total capacity    
		& 11.94 MVA 
		& 11.19 MVA    \\
		\bottomrule 
		\bottomrule 
	\end{tabular}
	%\vspace{-0.6cm}
\end{table}

\subsubsection{Converter capacity and cost savings}

The conventional ReP2H system requires two separate AC/DC converters for the electrolyzer and the ES, whereas the proposed ESEHE uses a shared AC/DC interface, eliminating the AC/DC stage dedicated to the ES. Considering that simultaneous full-power operation of the ES and electrolyzer occurs with very low probability, the optimal AC/DC converter capacity in the ESEHE can be smaller than the sum of two independent converters in conventional systems. A techno-economic evaluation is therefore conducted.

For both systems, 3\% converter losses are considered. In a MW-rated system, statistical analysis of annual PV data  (Fig. \ref{fig:11}) shows that power exceeds 5 MW only eight times (2\% of samples), with an average value of 5.07 MW. By balancing converter capacity cost against acceptable PV power clipping, simulation-based sizing results are presented in Table \ref{table:5}. Although the 5.22 MVA shared AC/DC interface slightly exceeds the 5.15 MVA electrolyzer converter in the conventional system, the total AC/DC and DC/DC converter capacity decreases from the original 5.97 MVA. This design prioritizes generation probability rather than cumulative peak ratings. By accepting minimal power clipping, the strategy reduces capital investment and lowers the overall hydrogen production cost.

With 1 MW-rated AC/DC and DC/DC converters priced at 404,000 CNY/MW and 443,000 CNY/MW \cite{en18010048, zeng2025optimal}, respectively, the conventional system requires 5.06 million CNY. In contrast, the proposed ESEHE requires 4.75 million CNY, saving 0.31 million CNY. {\color{black}Although distributed allocation of ES increases system topology complexity, the reduction in total converter capacity leads to measurable capital savings, and this trade-off merits further investigation.}

Note that optimal converter sizing depends on renewable penetration levels. For example, under higher renewable capacities \cite{10814660,yang2025achieving} (e.g., a 1:1.5 electrolysis-to-renewable ratio), the AC/DC converter capacity increases to extend full-load operation, and optimal sizing also requires further study.

\section{Conclusion}
\label{sec:Conclusion}

This paper presents the topology and control framework of an ESEHE. By decentralizing ES at electrolyzer DC buses, the proposed approach improves both robustness and economic performance in off-grid ReP2H systems. Key findings are:

1) The ESEHE enhances GFM capability, reducing frequency deviation and voltage fluctuation, outperforming conventional GFM HEs without ES.

2) Coordinated use of the electrolytic load, ES, and EDL effect maximizes transient support and maintains stable power transfer under $N$-1 contingencies. {\color{black}This enables reliable GFM support in harsh off-grid environments.}

3) Eliminating bidirectional AC/DC conversions compared with conventional systems improves energy efficiency by 0.23\% {\color{black}and reduces total investment cost on AC/DC and DC/DC converters by 6\%.}

4) The frequency-splitting strategy reduces ES degradation by 19.7\%, and SOC equalization further extends battery lifespan in multi-electrolyzer systems.

Future work will explore optimal converter and ES sizing in ESEHE-based ReP2H systems, balancing dynamic performance and cost, and optimizing degradation trade-offs between ES and electrolysis stacks. {\color{black}By integrating degradation, infrastructure investment, and maintenance cost modeling with techno-economic optimization, future work is expected to identify cost–performance optimal configurations, extend component lifespans, and further reduce system cost.}

\bibliographystyle{ieeetr}
%\bibliography{ref}

\end{document}